\documentstyle{amsppt}
\magnification=1200
\hoffset=-0.5pc
\nologo
\vsize=57.2truepc
\hsize=38.5truepc

\spaceskip=.5em plus.25em minus.20em

\define\Bobb{\Bbb}
\define\fra{\frak}

 \define\abramars{1}
 \define\almekump{2}
 \define\armcusgo{3}
 \define\armamonc{4}
     \define\atiy{5}
 \define\beilsche{6}
 \define\bkoucone{7}
 \define\brunvett{8}
 \define\canawein{9}
\define\chasetwo{10}
\define\chemltwo{11}
\define\cheveile{12}
\define\craifern{13}
\define\debarone{14}
\define\diracone{15}
\define\diractwo{16}
\define\dreznara{17}
   \define\fahim{18}
\define\fuchsste{19}
\define\gormacon{20}
\define\gotaythr{21}
\define\guistetw{22}
\define\hermaone{23}
 \define\herzone{24}
\define\higgmack{25}
\define\hinschon{26}
\define\hochsfiv{27}
\define\hochkoro{28}
\define\poiscoho{29}
 \define\souriau{30}
    \define\srni{31}
\define\claustha{32}
      \define\bv{33}
\define\extensta{34}
 \define\duality{35}
  \define\banach{36}
\define\oberwork{37}
 \define\twilled{38}
 \define\kaehler{39}
  \define\severi{40}
      \define\qr{41}
\define\illusboo{42}
\define\jacobfiv{43}
\define\kambtond{44}
\define\kasstora{45}
\define\kempfone{46}
\define\kemneone{47}
\define\kirilone{48}
\define\kirwaboo{49}
\define\kosmagri{50}
\define\kostaone{51}
\define\kostetwo{52}
\define\kralyvin{53}
\define\lermonsj{54}
   \define\lichn{55}
\define\mackbook{56}
    \define\mack{57}
\define\marswein{58}
\define\masspete{59}
\define\naramtwo{60}
\define\nelsoboo{61}
\define\palaione{62}
    \define\prad{63}
\define\praditwo{64}
\define\ramadthr{65}
\define\rinehone{66}
   \define\segal{67}
\define\sjamatwo{68}
\define\sjamafou{69}
\define\sjamlerm{70}
\define\sniabook{71}
\define\sniatone{72}
\define\sniawein{73}
\define\sourione{74}
\define\stashfiv{75}
\define\nteleman{76}
\define\ctelethr{77}
\define\weylbook{78}
\define\woodhous{79}

\topmatter
\title Lie-Rinehart algebras, descent, and quantization
\endtitle
\author Johannes Huebschmann
\endauthor
\affil 
Universit\'e des Sciences et Technologies
de Lille
\\
U. F. R. de Math\'ematiques
\\
CNRS-UMR 8524
\\
F-59 655 VILLENEUVE D'ASCQ, France
\\
Johannes.Huebschmann\@AGAT.UNIV-LILLE1.FR
\endaffil
\date{February 24, 2003}
\enddate
\abstract{A {\it Lie-Rinehart algebra\/} 
$(A,L)$ consists of a commutative algebra
$A$ and a Lie algebra $L$ with additional structure 
which generalizes the mutual structure of interaction
between the algebra of smooth functions and the Lie algebra
of smooth vector fields on a smooth manifold.
Lie-Rinehart algebras
provide the correct
categorical language
to solve the problem whether K\"ahler quantization commutes
with reduction which, in turn, may be seen as a descent problem.}
\endabstract

\address{\smallskip
\noindent
USTL, UFR de Math\'ematiques, 
CNRS-UMR 8524, F-59 655 Villeneuve d'Ascq C\'edex,
France
\newline\noindent
Johannes.Huebschmann\@agat.univ-lille1.fr}
\endaddress
\endkeywords

\endtopmatter
\document
\leftheadtext{Johannes Huebschmann}
\rightheadtext{Lie-Rinehart algebras, descent, and quantization}

\medskip
\noindent
{\bf Introduction}
\smallskip\noindent
The algebra of smooth functions $C^{\infty}(N)$
on a smooth manifold $N$ and its Lie algebra of smooth vector fields
$\roman{Vect}(N)$ have an interesting structure of interaction.
For reasons which will become apparent below,
we will refer to a pair $(A,L)$ which consists of a commutative algebra
$A$ and a Lie algebra $L$ with additional structure modeled on 
a pair of the kind $(C^{\infty}(N),\roman{Vect}(N))$ as a
{\it Lie-Rinehart algebra\/}.
In this article we will show that the notion of Lie-Rinehart algebra
 provides the correct
categorical language to solve a problem which we will describe
shortly. Lie-Rinehart algebras occur in other areas of mathematics
as well; an overview will be given in Section 1 below.
\smallskip
According to a philosophy going back to {\smc Dirac\/},
the correspondence between a classical theory and its 
quantum counterpart should be based on an analogy between their 
mathematical structures. In one direction, this correspondence,
albeit not well defined,
is 
referred to
as {\it quantization\/}.
Given a classical system with constraints which, in turn, determine
what is called the {\it reduced system\/},
the question arises whether 
quantization {\it descends\/} to the reduced system
in such a way that, once the 
unconstrained system has been successfully quantized,
imposing the symmetries on the quantized unconstrained system
is equivalent to quantizing  the reduced system.
This question goes back to the early days of quantum mechanics
and appears already in  {\smc Dirac's}
work on the electron and positron \cite\diracone, \cite\diractwo.
\smallskip
In the framework of K\"ahler quantization,
the problem may be phrased as a {\it descent problem\/} and, indeed,
under favorable circumstances which, essentially come down to requiring that
the unreduced and reduced spaces be both ordinary
quantizable K\"ahler manifolds, the problem has been known for long to have
a solution 
\cite\guistetw\ 
which, among other things, involves a version of what is
referred to as 
{\smc Kempf's\/} {\it descent lemma\/} \cite\kempfone\ 
in geometric invariant theory;
see e.~g. \cite\dreznara\ and Remark 8.6 below.
In the present article we will advertise the idea  that 
the concept of {\it Lie-Rinehart algebra\/}
provides the appropriate categorical language to solve the problem,
spelled out as a {\it descent problem\/}, 
under suitable more general circumstances so that,
in a sense,
reduction after quantization is then equivalent to
quantization after reduction;
here the term \lq\lq descent\rq\rq\ should, perhaps, not be taken in too narrow
a sense.
\smallskip
Given a classical system, its dynamical behaviour being encapsulated in
a Poisson bracket among the classical observables,
according to {\smc Dirac\/}'s idea of correspondence between the classical
and quantum system,
the Poisson bracket should then be the classical analogue of the quantum 
mechanical commutator. Thus, on the physics side,
the Poisson bracket is a crucial piece of structure. 
Mathematically, it is a crucial piece of structure as well;
in particular, when the classical phase space involves singularities,
these may be understood in terms of the Poisson structure.  
More precisely,
when the classical phase space carries a stratified symplectic structure,
the Poisson structure 
encapsulates the mutual positions of the symplectic structures
on the strata. 
See e.~g. \cite\srni, \cite\claustha, \cite\oberwork,
\cite\kaehler.
\smallskip
Up to now,
the available methods 
have been insufficient to attack the problem of
quantization of reduced observables, 
once the reduced phase space is no longer a smooth manifold;
we will refer to this situation as the {\it singular case\/}.
The singular case is the rule rather than the exception.
For example, simple classical mechanical systems and
the solution spaces of classical field theories
involve singularities;
see e.~g. \cite\armcusgo\ and the references there.
In the presence of singularities,
restricting quantization
to a smooth open dense part, the \lq\lq top stratum\rq\rq, leads
to a loss of information and in fact to
inconsistent results,
cf. Section 4 of \cite\qr.
To overcome
these difficulties on the classical level,
in \cite\kaehler,
we isolated a certain class
of \lq\lq K\"ahler spaces with singularities\rq\rq,
which we call {\it stratified K\"ahler spaces\/}.
On such a space,
the complex analytic structure alone is unsatisfactory
for issues related with
quantization because it overlooks the requisite Poisson structures.
In \cite\qr\ we developed the K\"ahler quantization scheme
over (complex analytic) stratified K\"ahler spaces.
A suitable notion of
{\it prequantization\/}, phrased in terms of {\it prequantum modules\/} 
introduced in \cite\souriau,
yields the requisite representation of the 
Poisson algebra; in particular, this representation satisfies the
Dirac condition. 
A suitably defined concept of
stratified K\"ahler polarization then takes care of the 
irreducibility problem, as  does an ordinary polarization
in the smooth case.
Over a stratified space,
the appropriate quantum phase space is what we call a
{\it costratified\/} Hilbert space; this is a system of Hilbert spaces,
one for each stratum, which arises from quantization
on the closure of that stratum, the stratification
provides linear maps between these Hilbert spaces
reversing the partial ordering among the strata,
and these linear maps are compatible with the quantizations. 
The main result obtained in \cite\qr\ says
that, for a positive K\"ahler manifold
with a hamiltonian action of a compact Lie group,
when suitable additional conditions are imposed, reduction after 
quantization coincides with quantization after reduction in the sense that not 
only the reduced and unreduced quantum phase spaces correspond but the 
{\it invariant unreduced and reduced quantum observables as well\/}.
Examples abound; one such class of examples, 
involving holomorphic nilpotent orbits and in particular angular momentum
zero spaces,
has been  treated in \cite\qr.
A particular case thereof will be reproduced in Section 6 below,
for the sake of illustration.
\smallskip
A stratified polarization, 
see Section 5 below for details, is defined in terms
of an appropriate Lie-Rinehart algebra which,
for any Poisson algebra, serves as a replacement for 
the tangent bundle of a smooth symplectic manifold.
The question whether quantization commutes with reduction
includes the question whether what is behind the phrase
\lq\lq in terms of an appropriate Lie-Rinehart algebra\rq\rq\ 
descends to the reduced level. This hints at
interpreting this question as
a descent problem.
\smallskip
To our knowledge, the idea of Lie-Rinehart algebra 
was first used by {\smc Jacobson\/} in \cite \jacobfiv\ 
(without being explicitly identified as a structure in its own)
to study certain
 field extensions.
Thereafter this idea occurred in other areas including 
differential geometry and differential Galois theory.
More details will be given below.
\smallskip
I am indebted  to the organizers of the meeting for having given me
the chance to illustrate an application of  Lie-Rinehart algebras
to a problem phrased in a language  entirely different from that of
Lie-Rinehart algebras.
Perhaps one can build a general Galois theory 
including ordinary Galois theory, differential Galois
theory, and ordinary principal bundles,
in which Lie-Rinehart algebras
appear as certain objects which capture infinitesimal symmetries.

\medskip\noindent{\bf 1. Lie-Rinehart algebras}
\smallskip\noindent
Let $R$ be a commutative ring with 1 taken as ground ring which,
for the moment,
may be arbitrary.
For a commutative $R$-algebra $A$,
we denote by $\roman{Der}(A)$ 
the $R$-Lie algebra of
derivations of $A$,
with its standard Lie algebra structure.
An $(R,A)$-{\it Lie algebra\/}
\cite\rinehone\ 
is a Lie algebra $L$ over $R$ which acts
on (the left of) $A$ 
(by derivations)
and is also an $A$-module 
satisfying suitable compatibility conditions
which generalize the usual properties
of the Lie algebra of vector fields on a smooth manifold
viewed as a module over its ring of functions;
these conditions read
$$
\align
[\alpha,a\beta]  &= \alpha(a)\beta + a [\alpha,\beta],
\\
(a\alpha)(b) &= a (\alpha(b)),
\endalign
$$ 
where
$a,b \in A$ and $\alpha, \beta \in L$.
When the emphasis is
on the pair $(A,L)$
with the mutual structure
of interaction between $A$ and $L$,
we refer to the pair $(A,L)$ as
a
{\it Lie-Rinehart} algebra.
Given  an arbitrary
commutative algebra $A$ over $R$,
an obvious example of a 
Lie-Rinehart algebra
is the pair $(A,\roman{Der}(A))$,
with the obvious action of $\roman{Der}(A)$ on $A$ 
and obvious $A$-module structure on $\roman{Der}(A)$.
There is an obvious notion of morphism of
Lie-Rinehart algebras and, with this notion of
morphism, Lie-Rinehart algebras constitute a category.
More details
may be found in {\smc Rinehart}~\cite\rinehone\ 
and in our papers \cite\poiscoho\ and \cite\souriau.
\smallskip

We will now briefly spell out
some of the salient features of Lie-Rinehart algebras.
Given 
an $(R,A)$-Lie algebra $L$, its 
{\it universal algebra\/}
${(U(A,L),\iota_L,\iota_A)}$
is an $R$-algebra $U(A,L)$ together with a morphism
${
\iota_A
\colon 
A
\longrightarrow
U(A,L)
}$
of $R$-algebras
and
a morphism
${
\iota_L
\colon 
L
\longrightarrow
U(A,L)
}$
of Lie algebras over $R$
having the properties
$$
\iota_A(a)\iota_L(\alpha) 
= \iota_L(a\,\alpha),\quad
\iota_L(\alpha)\iota_A(a) - \iota_A(a)\iota_L(\alpha) 
= \iota_A(\alpha(a)),
$$
and
${(U(A,L),\iota_L,\iota_A)}$
is {\it universal\/} among triples
${(B,\phi_L,\phi_A)}$
having these properties.
For example,
when
 $A$ is the algebra of smooth functions on a smooth
manifold $N$ and  $L$  the Lie algebra of 
smooth vector fields
on $N$, then $U(A,L)$ is the {\it algebra of 
(globally defined)
differential operators
on\/} $N$.
An explicit construction for the 
$R$-algebra ${U(A,L)}$ is given 
in {\smc Rinehart} \cite\rinehone.
See our paper \cite\poiscoho\ for an alternate construction
which employs the {\smc Massey-Peterson\/} \cite\masspete\  algebra.
\smallskip
The universal algebra $U(A,L)$ admits an
obvious filtered algebra structure
$U_{-1} \subseteq U_{0}\subseteq U_{1}\subseteq \dots $,
cf. {\smc Rinehart} \cite\rinehone,
where
$U_{-1}(A,L) = 0$
and where, for $p \geq 0$,
$U_p(A,L)$ is the left $A$-submodule of
$U(A,L)$ generated by products of at most $p$ elements of
the image  ${\overline L}$ of $L$ in $U(A,L)$,
and the associated graded object $E^0(U(A,L))$ 
inherits a commutative graded $A$-algebra structure.
The Poincar\'e-Birkhoff-Witt
Theorem for $U(A,L)$
then takes the following form where $S_A[L]$ 
denotes the symmetric $A$-algebra on $L$,
cf. (3.1) of {\smc Rinehart} \cite\rinehone.

\proclaim{Theorem 1.1 [Rinehart]}
For an $(R,A)$-Lie algebra $L$ 
which 
is projective
as an $A$-module,
 the canonical
$A$-epimorphism
$S_A[L]
\longrightarrow
E^0(U(A,L))
$
is an isomorphism of $A$-algebras.
\endproclaim

Consequently,
for an $(R,A)$-Lie algebra $L$ 
which 
is projective
as an $A$-module,
 the morphism
$\iota_L
\colon
L
\longrightarrow
U(A,L)
$ is injective.

\smallskip
The  construction of the ordinary Koszul complex
computing Lie algebra cohomology
carries over
as well:
Let $\Lambda_A(sL)$ be the exterior Hopf algebra over $A$
on the suspension $sL$ of $L$,
where \lq\lq suspension\rq\rq\ 
means that $sL$ is $L$ except that its elements are
regraded up by 1.
{\smc Rinehart} \cite\rinehone\ has proved that
the ordinary Chevalley-Eilenberg operator
induces an ${U(A,L)}$-linear operator $d$ on
$U(A,L) \otimes _A \Lambda_A(sL)$
(this is not obvious since $L$ is not an ordinary $A$-Lie algebra
unless $L$ acts trivially on $A$)
having square zero.
We will refer to
$$
K(A,L) = (U(A,L) \otimes _A \Lambda_A(sL),d)
\tag1.2
$$
as the
{\it Rinehart complex\/} for $(A,L)$.
It is manifest that the Rinehart complex is functorial
in $(A,L)$. 
Moreover, as a graded $A$-module, the resulting chain complex
$\roman{Hom}_{U(A,L)}(K(A,L),A)$ underlies the $A$-algebra 
$\roman{Alt}_A(L,A)$
of $A$-multilinear
functions on $L$ but, beware, the differential is linear only over the ground
ring $R$ and turns $\roman{Alt}_A(L,A)$ into a 
differential graded cocommutative algebra over $R$;
we will refer to this differential graded $R$-algebra as the
{\it Rinehart algebra\/} of $(A,L)$.
Rinehart also noticed that,
when $L$ is projective or free as a left $A$-module,
$K(A,L)$ is a
projective 
or free
resolution of
$A$ in the category of left $U(A,L)$-modules
according as $L$ is a projective or free
left $A$-module;
details may be found in \cite\rinehone.
In particular, the Rinehart algebra
$(\roman{Alt}_A(L,A),d)$ then computes the Ext-algebra
$\roman{Ext}^*_{U(A,L)}(A,A)$.
\smallskip
Rinehart also noticed that, when $A$ is the algebra of smooth functions on a 
smooth 
manifold $N$ and $L$ the Lie algebra of smooth vector fields
on $N$, then $(\roman{Alt}_A(L,A),d)(=\roman{Hom}_{U(A,L)}(K(A,L),A))$ 
is the ordinary {\it de Rham complex of\/} $N$
whence, as an algebra, the de Rham cohomology of $N$ amounts to the Ext-algebra
$\roman{Ext}^*_{U(A,L)}(A,A)$ over
the algebra $U(A,L)$ of differential operators on $N$.
Likewise, for a Lie algebra $L$ over $R$
acting trivially on $R$,
$K(R,L)$
is the ordinary Koszul complex;  
in particular, when $L$ is projective as an $R$-module,
${K(R,L)}$ is the ordinary Koszul resolution of the ground ring $R$.
{\sl Thus the cohomology of Lie-Rinehart algebras comprises
de Rham- as well as Lie algebra cohomology.\/}
In particular, this offers a possible explanation 
why {\smc Chevalley\/}  and {\smc Eilenberg\/}
\cite\cheveile, when they first isolated Lie algebra
cohomology, derived their formulas by abstracting from
the de Rham operator of a smooth manifold.
Suitable graded versions of the cohomology of Lie-Rinehart algebras
comprise as well
Hodge cohomology and coherent sheaf cohomology of complex manifolds
\cite{\banach,\,\twilled}.

\smallskip
The classical differential geometry notions of connection,
curvature,  characteristic classes, etc.
may be developed for arbitrary Lie-Rinehart algebras
\cite\poiscoho, \cite\souriau, \cite\extensta,
and there are  notions of  duality
 for Lie-Rinehart algebras
generalizing Poincar\'e duality \cite\duality; 
the idea of duality  has been shown in
\cite\bv\ to cast new light on Gerstenhaber- and Batalin-Vilkovisky 
algebras. 
In 
a sense
these homological algebra interpretations of 
Batalin-Vilkovisky algebras push further Rinehart's observations 
related with the interpretation of de Rham cohomology  
as certain Ext-groups.
Graded versions of duality for Lie-Rinehart algebras
\cite\banach, \cite\twilled\ 
may be used to study e.~g. complex manifolds,
CR-structures, and the mirror conjecture.
\smallskip
Lie-Rinehart algebras were implicitly used 
already by {\smc Jacobson\/} \cite\jacobfiv\ 
and later by {\smc Hochschild\/} \cite\hochsfiv.
The idea of Lie-Rinehart algebra has been 
introduced by a very large number of authors, most
of whom independently proposed their own terminology. 
I am indebted to K. Mackenzie
for his help with compiling the following list
in chronological order:
   Pseudo-alg\`ebre de Lie:  Herz, 1953 \cite\herzone---actually, Herz seems
to be the first to describe  the structure  in a form which makes its
generality clear---;
   Lie d-ring: Palais, 1961 \cite\palaione;
   (R,C)-Lie algebra: Rinehart, 1963 \cite\rinehone;
   (R,C)-\'espace d'Elie Cartan r\'egulier et sans courbure: 
   de Barros, 1964 \cite\debarone;
   (R,C)-alg\`ebre de Lie: Bkouche, 1966 \cite\bkoucone;
   Lie algebra with an associated module structure: Hermann, 1967 \cite\hermaone; 
   Lie module: Nelson, 1967 \cite\nelsoboo;
   Pseudo-alg\`ebre de Lie: Pradines, 1967 \cite\praditwo;
   (A, C) system: Kostant and Sternberg, 1971 \cite\kostetwo;
   Sheaf of twisted Lie algebras: Kamber and Tondeur, 1971 \cite\kambtond;
   Alg\`ebre de Lie sur C/R: Illusie, 1972 \cite\illusboo;
   Lie algebra extension: N. Teleman, 1972 \cite\nteleman;
   Lie-Cartan pair: Kastler and Stora, 1985 \cite\kasstora; 
   Atiyah algebra: Beilinson and Schechtmann, 1988 \cite\beilsche;
   Lie-Rinehart algebra: Huebschmann, 1990 \cite\poiscoho; 
   Differential Lie algebra: Kosmann--Schwarzbach and 
   Magri, 1990 \cite\kosmagri.
   Hinich and Schechtman (1993) \cite\hinschon\  have used the term   
   Lie algebroid for the general algebraic concept.
In differential Galois theory, 
Lie-Rinehart algebras occur under the name
\lq\lq alg\`ebre diff\'erentielle\rq\rq\ in a paper by 
{\smc Fahim\/} \cite\fahim.
Lie-Rinehart algebras occur as well in {\smc Chase\/} \cite\chasetwo\ 
and  {\smc Stasheff\/} \cite\stashfiv.
We have chosen to use the terminology {\it Lie-Rinehart algebra\/}
since, as already pointed out, Rinehart subsumed the
cohomology of these objects under standard homological algebra
and established a Poincar\'e-Birkhoff-Witt theorem for 
them.
In differential geometry, $(R,A)$-Lie algebras arise as spaces of 
sections of Lie algebroids. These, in turn, 
were introduced in 1966 by {\smc Pradines\/} \cite\prad\
and, in that paper, Pradines raised the issue whether
Lie's third theorem holds for Lie algebroids in the sense that
any Lie algebroid integrates to a Lie groupoid.
{\smc Crainic\/} and {\smc Fernandes\/} \cite\craifern\ 
have recently given a solution of this problem in terms of suitably
defined obstructions.
See {\smc Mackenzie\/} \cite\mackbook\ for a complete account of Lie algebroids and
Lie groupoids, as well as 
{\smc Canas da Silva-Weinstein\/}
\cite\canawein\ and  {\smc Mackenzie\/} \cite\mack\
for more recent surveys on particular aspects.
The idea of Lie algebroid is lurking behind a construction in 
{\smc Fuchssteiner\/} \cite\fuchsste\ 
(see Remark 2.6 below).
A descent construction for Lie algebroids may be found in
{\smc Higgins-Mackenzie\/}
\cite\higgmack.
A general notion of morphism of Lie algebroids has been introduced
by {\smc Almeida-Kumpera\/}
\cite\almekump. This notion has been used by 
S. {\smc Chemla\/} \cite\chemltwo\ to study
notions of duality for Lie algebroids in complex algebraic geometry generalizing
Serre duality. 
Lie-Rinehart algebras are lurking as well behind the nowadays very
active research area of $D$-modules.
\smallskip\noindent
{\smc Remark} (out of context).
At the end of his paper \cite\rinehone, Rinehart  introduced an operator
on the Hochschild complex of a commutative algebra which, 
some 20 years later,
was reinvented by A. Connes in order to define cyclic cohomology.

\medskip\noindent{\bf 2. Poisson algebras}
\smallskip\noindent
For intelligibility, we recall briefly how
for an arbitrary Poisson algebra
an appropriate Lie-Rinehart algebra serves as a replacement for 
the tangent bundle of a smooth symplectic manifold.

\smallskip
Let $(A,\{\cdot,\cdot\})$ be a Poisson algebra, 
and let $D_A$ be the $A$-module
of formal differentials of $A$ the elements of which we write
as $du$, for $u \in A$.
For $u,v \in A$, the assignment to $(du,dv)$ of
$
\pi (du,dv) = \{u,v\}
$
yields an $A$-valued $A$-bilinear skew-symmetric 2-form 
$\pi= \pi_{\{\cdot,\cdot\}}$
on $D_A$,
the {\it Poisson\/} 2-{\it form\/} 
for 
$(A,\{\cdot,\cdot\})$.
Its adjoint 
$$
\pi^{\sharp} 
\colon D_A
@>>> 
\roman{Der}(A) = \roman{Hom}_A(D_A,A)
\tag2.1
$$
is a morphism of $A$-modules,
and the formula
$$
[a du,b dv] = 
a\{u,b\} dv + b\{a,v\} du + ab d\{u,v\}
\tag2.2
$$
yields a Lie bracket
$[\cdot,\cdot]$
on $D_A$, viewed as an $R$-module.
More details may be found in
\cite\poiscoho.
For the record we recall the following, 
established in \cite\poiscoho\ (3.8).

\proclaim{Proposition 2.3}
The $A$-module structure on $D_A$,
the bracket  $[\cdot,\cdot]$,
and the morphism
 $\pi^{\sharp}$
of $A$-modules 
turn the pair
 $(A,D_A)$
into a 
Lie-Rinehart algebra  
in such a way that
$\pi^{\sharp}$
is a morphism of Lie-Rinehart algebras.
\endproclaim

We write
$D_{\{\cdot,\cdot\}} = (D_A,[\cdot,\cdot],\pi^{\sharp})$.
The
2-form
$
\pi_{\{\cdot,\cdot\}},
$
which is defined for {\it every\/} Poisson algebra,
is plainly a 2-cocycle in the Rinehart algebra
$(\roman{Alt}_A(D_{\{\cdot,\cdot\}},A),d)$.
In \cite\poiscoho,
we defined the {\it Poisson cohomology\/} 
$\roman H^*_{\roman{Poisson}}(A,A)$
of the Poisson algebra
$(A,\{\cdot,\cdot\})$ 
to be the cohomology of this Rinehart algebra, that is,
$$
\roman H^*_{\roman{Poisson}}(A,A)
=\roman H^*\left(\roman{Alt}_A(D_{\{\cdot,\cdot\}},A),d\right).
$$
\smallskip
Henceforth we shall take as ground ring that of the reals
$\Bobb R$ or that of the complex numbers $\Bobb C$.
We shall consider spaces $N$ with an algebra
of continuous $\Bobb R$-valued or $\Bobb C$-valued functions, 
deliberately denoted by
$C^{\infty}(N,\Bobb R)$ or $C^{\infty}(N,\Bobb C)$ as appropriate,
or by just $C^{\infty}(N)$,
for example ordinary smooth manifolds and ordinary smooth functions;
such an algebra 
$C^{\infty}(N)$
will then be referred to as a {\it smooth\/} structure
on $N$,
and 
$C^{\infty}(N)$
will be
viewed 
as part of the structure of $N$.
A space may support
{\it different\/}
smooth structures, though.
Given a space $N$
with a smooth structure $C^{\infty}(N)$,
we shall write
$\Omega^1(N)$ for the space of formal differentials with
those differentials divided out that are zero at each point,
cf. \cite\kralyvin;
for example, over the real line with its ordinary smooth structure, 
the formal differentials
$d \sin x$ and $ \cos x dx$ do not coincide but the
formal differential
$d \sin x - \cos x dx$
is zero at each point.
At a point of $N$, 
the object
$\Omega^1(N)$
comes down to the ordinary 
space of differentials for the smooth structure on $N$; see Section 1.3 of 
our paper \cite\qr\ for details.
When 
$N$ is an ordinary smooth manifold,
$\Omega^1(N)$ amounts to the space of smooth sections
of the cotangent bundle.
For a general smooth space $N$ over the reals, when
$A= C^{\infty}(N,\Bobb R)$
is endowed with a Poisson structure,
the  formula (1.2) yields a Lie-bracket
$[\cdot,\cdot]$ on
the $A$-module $\Omega^1(N)$
and the 2-form $\pi_{\{\cdot,\cdot\}}$ 
is still defined on 
$\Omega^1(N)$; its 
adjoint  
then yields an $A$-linear map
$
\pi^{\sharp}
$ 
from
$\Omega^1(N)$
to
$\roman{Der}(A)$ in such a way
that
$([\cdot,\cdot],\pi^{\sharp})$
is an $(\Bobb R,A)$-Lie algebra structure
on $\Omega^1(N)$
and that
the adjoint  
$
\pi^{\sharp}
$ 
is a morphism of $(\Bobb R,A)$-Lie algebras.
Here the notation
$[\cdot,\cdot]$, $\pi^{\sharp}$, $\pi_{\{\cdot,\cdot\}}$
is abused somewhat.
The obvious projection map
from
$D_{C^{\infty}(N)}$ to
$\Omega^1(N)$
is plainly compatible with the  Lie-Rinehart structures.
The fact that $D_{C^{\infty}(N)}$ is \lq\lq bigger\rq\rq\  than
$\Omega^1(N)$
in the sense that the surjection from 
the former to the latter
has a non-trivial kernel
causes no problem here since
the $A$-dual of this surjection, that is, the induced map from
$\roman{Hom}(\Omega^1(N),A)$ to
$\roman{Hom}(D_{C^{\infty}(N)},A)$, is an isomorphism.
We shall write
$\Omega^1(N)_{\{\cdot,\cdot\}}=(\Omega^1(N),[\cdot,\cdot],\pi^{\sharp})$.
When $N$ is a smooth manifold in the  usual sense,
the range 
$\roman{Der}(A)$
of the  adjoint map
$\pi^{\sharp}$
from
$\Omega^1(N)$ to
$\roman{Der}(A)$
boils down to the space
$\roman{Vect}(N)$
of smooth vector fields on $N$.
In this case, 
the Poisson structure 
on $N$ is symplectic,
that is, arises from a (uniquely determined) symplectic structure
on $N$, if and only if $\pi^{\sharp}$,
which may now be written
as a morphism
of smooth vector bundles from $\roman T^*N$ to $\roman TN$,
is an isomorphism.
\smallskip
The
2-form
$
\pi_{\{\cdot,\cdot\}}
$
is  a 2-cocycle in the Rinehart algebra
$(\roman{Alt}_A(\Omega^1(N)_{\{\cdot,\cdot\}},A),d)$ and
the canonical map
of differential graded algebras from
$(\roman{Alt}_A(\Omega^1(N)_{\{\cdot,\cdot\}},A),d)$
to $(\roman{Alt}_A(D_{\{\cdot,\cdot\}},A),d)$
is an isomorphism.
In particular, we may take the cohomology of the Rinehart algebra
$(\roman{Alt}_A(\Omega^1(N)_{\{\cdot,\cdot\}},A),d)$
as the definition of the Poisson cohomology
$\roman H^*_{\roman{Poisson}}(A,A)$
of the Poisson algebra
$(A,\{\cdot,\cdot\})$ 
as well.
When $N$ is an ordinary smooth manifold, its algebra of ordinary
smooth functions being endowed with a Poisson structure,
this notion of Poisson cohomology comes down to 
that introduced by {\smc Lichnerowicz\/} \cite\lichn.
For a general 
Poisson algebra $A$,
the
2-form
$
\pi_{\{\cdot,\cdot\}},
$
be it defined on $\Omega^1(N)_{\{\cdot,\cdot\}}$
for the case where $A$ is the structure algebra  
$C^{\infty}(N)$ of a space $N$ or on 
$D_{\{\cdot,\cdot\}}$ for an arbitrary Poisson algebra,
generalizes the symplectic form of a symplectic manifold;
see Section 3 of~\cite\poiscoho\  for details.
Suffice it to make the following observation, relevant for quantization:
Consider a space $N$ with a smooth structure
$C^{\infty}(N)$ which, in turn, is endowed with a Poisson bracket
$\{\cdot,\cdot\}$.
The
Poisson 2-form $\pi_{\{\cdot,\cdot\}}$
determines an {\it extension\/}
of Lie-Rinehart algebras which is central as a Lie algebra extension.
For technical reasons it is more convenient to take here the extension
$$
0 
@>>> 
A
@>>> 
\overline L_{\{\cdot,\cdot\}}
@>>> 
\Omega^1(N)_{\{\cdot,\cdot\}}
@>>> 
0
\tag2.4
$$
which corresponds to the negative of the Poisson 2-form.
Here, as $A$-modules,
$\overline L_{\{\cdot,\cdot\}} = A \oplus \Omega^1(N)_{\{\cdot,\cdot\}}$,
and the Lie bracket on $\overline L_{\{\cdot,\cdot\}}$ 
is given by
$$
[(a,du),(b,dv)] =
\left(
\{u,b\}+ \{a,v\} - \{u,v\}, d\{u,v\}
\right) ,\quad
a,b,u,v \in A.
\tag2.5
$$
Here we have written
\lq\lq $\overline L$\rq\rq\ 
rather than simply $L$
to indicate that
the extension (2.4) represents the negative of the class of
$\pi_{\{\cdot,\cdot\}}$
in the second  cohomology group
$\roman H^2(\roman{Alt}_A(D_{\{\cdot,\cdot\}},A),d)$ of the corresponding
Rinehart algebra,
cf. \cite\poiscoho, and the notation $du$, $dv$ etc. is abused somewhat.
Now,
any principal circle bundle admits as its infinitesimal object
an {\it Atiyah\/}-sequence
\cite\atiy\ 
whose spaces of sections
constitute
a central extension of Lie-Rinehart algebras;
see \cite\mackbook\ for a complete account of Atiyah-sequences
and \cite\extensta\ for a theory of characteristic classes for
 extensions of general Lie-Rinehart algebras.
When the Poisson structure is an ordinary smooth
symplectic Poisson structure whose symplectic form
represents an integral cohomology class, 
the Lie-Rinehart algebra extension
(2.4) comes down to the space of sections of the Atiyah-sequence
of the principal circle bundle classified by that cohomology class. 

\smallskip\noindent
{\smc Remark 2.6.}
For the special case where
$N$ is an ordinary smooth manifold, $C^{\infty}(N)$
its algebra of ordinary smooth functions,
and where $\{\cdot,\cdot\}$ is a Poisson structure on
$C^{\infty}(N)$,
the Lie-Rinehart structure
on the pair $(C^{\infty}(N),\Omega^1(N))$
(where $\Omega^1(N)$ amounts to the space of ordinary smooth 1-forms on $N$)
was discovered by a number of authors during the 80's most of whom phrased
the structure in terms of the corresponding Lie algebroid
structure on the cotangent bundle of $N$;
some historical comments may be found in \cite\poiscoho.
The first reference 
I am aware of where versions of the Lie algebroid bracket
and of the anchor map
may be found  is \cite\fuchsste; in that paper,
the notion of \lq\lq implectic operator\rq\rq\
is introduced---this is the operator nowadays referred to
as {\it Poisson tensor\/}---and the Lie bracket and anchor map
are the formula (2) and morphism written as $\Omega_{\phi}$,
respectively, in that paper.
{\sl The construction in terms of formal differentials
carried out in\/} \cite\poiscoho\
(and reproduced above)---as opposed to the
Lie algebroid construction---is {\sl
more general, though, since it refers to an arbitrary
Poisson structure, not necessarily one which is defined on an algebra
of smooth functions on an ordinary smooth manifold\/}.
In fact, the aim of the present article is to demonstrate the 
significance of this more general construction which works as well
for Poisson algebras of continuous functions defined on
 spaces with {\it singularities\/}
where among other things it yields a tool to relate the Poisson
structures on the strata of a stratified symplectic space; 
suitably translated
into the language of sheaves, it also works over not necessarily
non-singular varieties.

\medskip\noindent{\bf 3. Quantization}
\smallskip\noindent
According to {\smc Dirac}~\cite\diracone,~\cite\diractwo,
a {\it quantization\/} of a classical system 
described by a real Poisson algebra
$(A,\{\cdot,\cdot\})$ 
is a representation
$a \mapsto \hat a$
of a certain Lie subalgebra
$B$ of $A$,
$A$ and $B$ 
being viewed merely as Lie algebras,
by symmetric or, whenever possible, self-adjoint, operators
$\hat a$
on a Hilbert space
$\Cal H$ such that
(i) the Dirac condition
$$
i\,[\hat a,\hat b] = \widehat {\{a,b\}}
$$
holds;
that (ii)
for a  constant $c$, the operator $\hat c$ is given by
$
\hat c = c\, \roman{Id};
$
and that (iii)
the representation is irreducible.
Here the factor $i$ in the Dirac condition is forced by
the interpretation of quantum mechanics: Observables
are to be represented by symmetric (or self-adjoint) operators
but the ordinary commutator of two symmetric operators
is not symmetric.
The second requirement rules out the adjoint representation, and
the irreducibility condition is forced by the requiremend
that phase transitions be possible between two different states.
See e.~g. {\smc Sniatycki\/} \cite\sniabook\ or
{\smc Woodhouse}~\cite\woodhous.
Also it is known that for $B=A$ the problem has no solution
whence the requirement that only a sub Lie algebra of $A$ be 
represented.
The physical constant
$\hbar$ 
is here absorbed
in the Poisson structure.
It has become common to refer to a procedure
furnishing a representation that satisfies only (i) and (ii) above
as {\it prequantization\/}.
\smallskip
Under suitable circumstances,
over a smooth symplectic manifold,
the {\it geometric quantization scheme\/}
developed by
{\smc Kirilliov}~\cite\kirilone,
{\smc Kostant}~\cite\kostaone,
{\smc Souriau}~\cite\sourione,
and {\smc I. Segal}~\cite\segal,
furnishes a quantization; see \cite\sniabook\  or \cite\woodhous\ 
for complete accounts.
We confine ourselves with the remark that geometric quantization
proceeds in two steps. The first step, prequantization,
yields a representation of the Lie algebra underlying the whole Poisson
algebra which satisfies (i) and (ii) but such a representation is not
irreducible;
the second step involves a choice of {\it polarization\/}
to force the irreducibility condition.
In particular, a {\it K\"ahler polarization\/} leads to what is called
{\it K\"ahler quantization\/}. The existence of a K\"ahler polarization
entails that the underlying manifold carries an ordinary K\"ahler structure.
In the singular case,
the ordinary geometric quantization scheme is
no longer available, though.
In the rest of the paper we shall describe how,
under certain favorable circumstances, the difficulties
in the singular case can be overcome in the framework of K\"ahler
quantization.
An  observation crucial in the singular case 
is that the notion of polarization can be given
a meaning  by means of appropriately defined
Lie-Rinehart algebras.
Before going into details, we will briefly explain one of the origins of
singularities.

\medskip\noindent{\bf 4. Symmetries}
\smallskip\noindent
Recall that a {\it symplectic\/} manifold is a smooth manifold $N$ 
together
with a closed non-degenerate 2-form $\sigma$.
Given a function  $f$, the identity $\sigma(X_f,\cdot ) = df$ then associates a uniquely
determined vector field $X_f$ to $f$, the {\it Hamiltonian vector field\/}
of $f$ and, given two functions
$f$ and $h$,
their  
{\it Poisson bracket\/} $\{f,h\}$
is defined by $\{f,h\} = X_f h$. This
yields a Poisson bracket $\{\cdot,\cdot\}$ on the  algebra $C^{\infty} (N)$
of ordinary smooth functions on $N$, referred to as a {\it symplectic\/}
Poisson bracket.
\smallskip
Given a Lie group $G$,
a {\it hamiltonian\/} $G$-space is a smooth symplectic $G$-manifold 
$(N,\sigma)$
together with a smooth $G$-equivariant map $\mu$ from $N$ to the dual
$\fra g^*$ of the Lie algebra $\fra g$ of $G$ satisfying the formula
$$
\sigma(X_N,\cdot) = X \circ d \mu
\tag4.1
$$
for every $X \in \fra g$; here $X_N$ denotes the vector field on $N$
induced by $X \in \fra g$ via the $G$-action,
and $X$ is viewed as a linear form on $\fra g^*$.
The map $\mu$ is called a {\it momentum mapping\/}
(or {\it moment map\/}).
We  recall that
(4.1) says that,
given $X \in \fra g$,
the vector field $X_N$ is the hamiltonian vector field for 
the smooth function
$X\circ \mu$ on $N$.
See e.~g. \cite{\abramars} for details.
Given a hamiltonian $G$-space $(N,\sigma,\mu)$,
the space $N_{\roman{red}} = \mu^{-1}(0)\big / G$
is called the {\it reduced space\/}.
When $G$ is not compact, this space may have bad properties;
for example, it is not even a Hausdorff space 
when there  are non-closed $G$-orbits in $N$.

\smallskip
Let
$C^{\infty}(N_{\roman{red}}) = (C^{\infty}(N))^G)\big/I^G$,
where
$I^G$
refers to the ideal (in the algebra
$(C^{\infty}(N))^G$)
of smooth $G$-invariant functions on $N$
which vanish on the zero locus $\mu^{-1}(0)$;
this is a smooth structure on the reduced space $N_{\roman{red}}$.
As observed by {\smc Arms-Cushman-Gotay\/} \cite\armcusgo, 
the Noether Theorem implies that the symplectic Poisson structure
on $C^{\infty}(N)$
descends to a Poisson structure 
$\{\ ,\ \}_{\roman{red}}$
on
$C^{\infty}(N_{\roman{red}})$, and
{\smc Sjamaar-Lerman\/} \cite\sjamlerm\ have shown that, when $G$
is compact and when the momentum mapping is proper,
the orbit type decomposition of
$N_{\roman{red}}$  is a stratification in the sense of 
{\smc Goresky-MacPherson\/}
\cite\gormacon.
The idea that the orbit type decomposition is a stratification
(in a somewhat weaker sense) may be found already in \cite\armamonc.
For intelligibility, we recall some of the requisite
technical details.
\smallskip
A decomposition of a space $Y$ into 
pieces which are
smooth manifolds
such that these pieces fit together in a certain precise way
is called a {\it stratification\/}
\cite\gormacon.
More precisely:
Let $Y$ be a Hausdorff paracompact topological space
and let $\Cal I$ be a partially ordered set
with order relation denoted by $\leq$.
An $\Cal I$-{\it decomposition \/} of $Y$ is a locally finite collection
of disjoint locally closed manifolds $S_i \subseteq Y$
called {\it pieces\/} 
(recall that a collection $\Cal A$ of subsets of $Y$ is said to be
{\it locally finite\/} provided every $x \in Y$ has a neighborhood
$U_x$ in $Y$ such that $U_x \cap A \ne \emptyset$
for at most finitely many $A$ in $\Cal A$)
such that the following hold:
\smallskip
$Y = \cup S_i\ (i \in \Cal I)$,
\smallskip
$S_i\cap \overline S_j \ne \emptyset \Longleftrightarrow S_i \subseteq 
\overline S_j \Longleftrightarrow i \leq j\ (i,j \in \Cal I)$.
\smallskip\noindent
The space $Y$ is then called a {\it decomposed\/} space.
A decomposed space $Y$ is said to be a {\it stratified space\/}
if the pieces of $Y$, called {\it strata\/},
satisfy
the following condition:
Given a point $x$ in a piece $S$ there is an open neighborhood $U$
of $x$ in $Y$,
an open ball $B$ around $x$ in $S$, a stratified space $\Lambda$, called the 
{\it link\/} of $x$, and a decomposition preserving homeomorphism
from $B \times C^{\circ}(\Lambda)$ onto $U$.
Here $C^{\circ}(\Lambda)$
refers to the open cone on $\Lambda$ and, as a stratified space,
$\Lambda$ is less complicated than $C^{\circ}(\Lambda)$
whence the definition is not circular;
the idea of complication is here made precise by means of the notion
of {\it depth\/}.
\smallskip
A {\it stratified symplectic space\/}
\cite\sjamlerm\ 
is a stratified space $Y$
together with a 
Poisson algebra $(C^{\infty}(Y),\{\ ,\ \})$
of continuous functions on $Y$ which, on each piece of the decomposition,
restricts to an ordinary smooth
symplectic Poisson structure;
in particular,
$C^{\infty}(Y)$ is a smooth structure on $Y$.

\smallskip\noindent
{\smc Example 4.2.}
On the ordinary plane, with coordinates $x_1,x_2$,
consider
the algebra $A$ of smooth functions in
the coordinate functions
$x_1$, $x_2$ together with, 
which is {\it crucial\/} here, an additional function
$r$ which is the radius function, subject to the relation 
$x_1^2 + x_2^2 = r^2$.
Notice that $r$ is {\it not\/} a smooth function in the usual sense
whence the algebra $A$ 
is
strictly larger than that of ordinary smooth functions on the plane.
The  Poisson structure $\{\cdot,\cdot\}$ on $A$ given by
the formulas
$$
\{x_1,x_2\} = 2r,
\quad
\{x_1,r\} = 2x_2,
\quad
\{x_2,r\} = -2x_1
\tag4.2.1
$$
turns the plane into a stratified symplectic space.
Geometrically, the plane is taken here as a half cone,
the algebra $A$ being that of {\it Whitney\/}-smooth functions
on the half cone, with reference to the embedding
into 3-space;
there are two strata, the vertex of the half cone and the complement thereof.
On the complement of the vertex, which is a punctured plane,
the Poisson structure is symplectic.
In physics, the Poisson algebra $(A,\{\cdot,\cdot\})$ 
arises, for $n \geq 2$, as the reduced Poisson algebra
of a single particle in  $\Bobb R^n$
with $\roman O(n,\Bobb R)$-symmetry and angular momentum zero.
For $n=1$, the example still makes sense: the symmetry group is then just
a copy of $\Bobb Z/2$, and the angular momentum is zero.

\smallskip
Given a hamiltonian $G$-space $(N,\sigma,\mu)$ with $G$ compact,
in view of an observation in \cite\sjamlerm,
the Arms-Cushman-Gotay construction turns
$(N_{\roman{red}}, C^{\infty}(N_{\roman{red}}),\{\ ,\ \}_{\roman{red}})$
(more precisely: each connected component of
$N_{\roman{red}}$ in case the momentum mapping is not proper)
into a stratified symplectic space.
When $N_{\roman{red}}$ is smooth, i.~e. has a single stratum,
this space
is just a smooth symplectic manifold,
the ordinary {\smc Marsden-Weinstein\/} reduced space \cite\marswein.

\beginsection 5. Stratified complex polarizations

Within the ordinary geometric quantization scheme,
the irreducibility requirement is taken care of by means of
a {\it polarization\/}.
In particular, a {\it complex polarization\/} 
for an ordinary symplectic manifold $N$
is 
an integrable Lagrangian distribution
$F \subseteq \roman T^{\Bobb C}N$ 
of the complexified tangent bundle
$\roman T^{\Bobb C}N$
\cite\woodhous;
under the identification of
$\roman T^{\Bobb C}N$ with its (complex) dual
coming from the symplectic structure,
a complex polarization $F$ then
corresponds to a certain uniquely defined
$(\Bobb C, C^{\infty}(N,\Bobb C))$-Lie subalgebra
$P$ of
$\Omega^1(N,\Bobb C)_{\{\cdot,\cdot\}}$.
\smallskip
Given a stratified symplectic space $N$,
we refer to a  $(\Bobb C,C^{\infty}(X,\Bobb C))$-Lie subalgebra
$P$ of $\Omega^1(X,\Bobb C)_{\{\cdot,\cdot\}}$
as
a {\it stratified complex polarization\/}
for $N$ if, 
for every 
stratum $Y$,
under the identification of
$\roman T^{\Bobb C}Y$ with its (complex) dual
coming from the symplectic structure on that stratum,
the $(\Bobb C, C^{\infty}(Y,\Bobb C))$-Lie subalgebra
$P_Y$ of
$\Omega^1(Y,\Bobb C)_{\{\cdot,\cdot\}}$
generated by the restriction of $P$ to $Y$
is identified with 
the space of sections of
an ordinary complex polarization.
A stratified complex polarization is, then, a
{\it K\"ahler\/} polarization
provided on any stratum it comes from an ordinary (not necessarily positive)
K\"ahler polarization.   
We say that a stratified K\"ahler polarization is 
{\it complex analytic\/} provided it is induced from a 
complex analytic structure
on $N$, 
and we define a
complex analytic stratified K\"ahler structure
to be a {\it normal K\"ahler structure\/}
provided the complex analytic structure is normal.
A normal K\"ahler structure is {\it positive\/}
provided it is positive on each stratum.
See Section 2 of
\cite\kaehler\ 
for more details.
\smallskip
Let $G$ be a compact Lie group, 
denote its complex form by $G^{\Bobb C}$, and 
recall the following, cf. Proposition 4.2 of \cite\kaehler.

\proclaim{Proposition 5.1}
Given a positive K\"ahler manifold $N$
with a holomorphic $G^{\Bobb C}$-action 
whose restriction to 
$G$
preserves the K\"ahler structure
and
is  hamiltonian, 
the K\"ahler polarization $F$
induces a positive normal
(complex analytic stratified) K\"ahler polarization $P^{\roman{red}}$ 
on the reduced space $N^{\roman{red}}$, the latter being endowed
with its stratified symplectic Poisson algebra
$(C^{\infty}(N^{\roman{red}}),\{\cdot,\cdot\}^{\roman{red}})$.
\endproclaim

Under these circumstances, the underlying complex analytic structure
of $N^{\roman{red}}$ is that of a
geometric invariant theory quotient;
the existence thereof may be found in
\cite\kemneone\ and \cite\kirwaboo.
The existence problem of this complex analytic structure may be
seen as one of descent.

\medskip\noindent {\bf 6. Examples}
\smallskip\noindent
We will now illustrate the notions introduced so far
by means of a number of examples. The interested reader
will find more details in \cite\kaehler.
\smallskip\noindent
{\smc Example 6.1.} For $\ell \geq 1$,
consider the constrained system
of $\ell$ particles 
in $\Bobb R^s$ 
with total angular momentum zero.
Its unreduced
phase space
$N$ is a product $(\roman T^*\Bobb R^s)^{\ell}$
of $\ell$ copies 
of $\roman T^*\Bobb R^s$,
and we write the points of $N$ in the form
$(\bold q_1,\bold p_1,\dots, \bold q_\ell,\bold p_\ell)$.
Let $H=\roman O(s,\Bobb R)$. With reference to the obvious $H$-symmetry,
the momentum mapping of this system has the form
$$
\mu \colon N @>>> \fra h^*,\quad
\mu(\bold q_1,\bold p_1,\dots, \bold q_\ell,\bold p_\ell) 
= \bold q_1 \wedge \bold p_1 +
\dots + \bold q_\ell \wedge \bold p_\ell,
$$
where the Lie algebra $\fra h = \fra {so}(s,\Bobb R)$ is identified
with its dual in the standard fashion.
To elucidate the reduced space,
observe that
the assignment to
$(\bold q, \bold p) =(\bold q_1,\bold p_1,\dots, \bold q_\ell,\bold p_\ell)$
of the real symmetric
$(2\ell \times 2\ell)$-matrix
$
\xi(\bold q,\bold p)=
\left[
\matrix 
\bold q_j \bold q_k
&
\bold q_j \bold p_k
\\
\bold p_j \bold q_k
&
\bold p_j \bold p_k 
\endmatrix
\right]_{1 \leq j, k \leq \ell}
$
yields a real algebraic map
$
\xi \colon
N
@>>>
\roman S_{\Bobb R}^2[\Bobb R^{2\ell}]
$
from $N$ to the 
real vector space 
$\roman S_{\Bobb R}^2[\Bobb R^{2\ell}]$
of real symmetric
$(2\ell \times 2\ell)$-matrices
which passes
to an embedding
of
the  
reduced space
$N^{\roman{red}}=\mu^{-1}(0) \big / H$
into 
$\roman S_{\Bobb R}^2[\Bobb R^{2\ell}]$, in fact,
realizes
$N^{\roman{red}}$ as a semi-algebraic set
in $\roman S_{\Bobb R}^2[\Bobb R^{2\ell}]$.
Let $J$ be the standard complex structure on
$\Bobb R^{2\ell}$.
Now, on the one hand,
the association
$S \mapsto JS$
identifies
$\roman S_{\Bobb R}^2[\Bobb R^{2\ell}]$
with 
$\fra{sp}(\ell,\Bobb R)$
in an $\roman {Sp}(\ell,\Bobb R)$-equivariant
fashion
(with reference to the obvious actions) 
and hence identifies
$N^{\roman{red}}$ with a subset
of $\fra{sp}(\ell,\Bobb R)$
which,
as observed in
\cite\lermonsj, 
is  the closure of a certain nilpotent orbit which
has been identified as a {\it holomorphic nilpotent orbit\/}
in \cite\kaehler.
The Killing form transforms
the Lie-Poisson structure on
$\fra{sp}(\ell,\Bobb R)^*$
to a Poisson structure
on
$\fra{sp}(\ell,\Bobb R)$
which, restricted to
$N^{\roman{red}}$,
yields a stratified symplectic structure.
Another observation in \cite\lermonsj\ 
entails that
this stratified symplectic structure
coincides with the
Sjamaar-Lerman
stratified symplectic structure
\cite\sjamlerm\ 
mentioned earlier
arising by symplectic reduction
from
$N=(\roman T^*\Bobb R^s)^{\ell}$. 
We mention in passing that,
$\fra{sp}(\ell,\Bobb R)$ being identified with its dual
by means of an appropriate positive multiple of the Killing form,
as well as with
$\roman S_{\Bobb R}^2[\Bobb R^{2\ell}]$,
the map $\xi$ is essentially the momentum mapping for the obvious
$\roman {Sp}(\ell,\Bobb R)$-action on $N$.

On the other hand,
the choice of $J$ determines a maximal compact subalgebra
of $\fra{sp}(\ell,\Bobb R)$ which is just a copy of
$\fra u(\ell)$ and, furthermore, a Cartan decomposition
$\fra{sp}(\ell,\Bobb R)= \fra u(\ell) \oplus \fra p$.
Now matrix multiplication by $J$ from the left induces a complex
structure on $\fra p$ and, with this structure,
as a complex vector space, $\fra p$ amounts to 
the complex symmetric square
$\roman S_{\Bobb C}^2[\Bobb C^{\ell}]$
on $\Bobb C^{\ell}$.
In particular, orthogonal projection to $\fra p$ induces
a linear surjection
of real vector spaces
from $\roman S_{\Bobb R}^2[\Bobb R^{2\ell}]$
to
$\roman S_{\Bobb C}^2[\Bobb C^{\ell}]$,
uniquely determined by $J$; it is
given by
the assignment to a real symmetric
$(2\ell \times 2\ell)$-matrix
of the corresponding
complex
symmetric
$(\ell \times \ell)$-matrix
with respect to the standard complex structure
$J$ on
$\Bobb R^{2\ell}$.
This projection,
restricted to
$N^{\roman{red}}$,
is injective and 
yields a complex analytic structure on 
$N^{\roman{red}}$.
The two structures are compatible
and yield a normal (complex analytic stratified)
K\"ahler structure on $N^{\roman{red}}$; see \cite\kaehler\ for details.
We will describe the requisite (complex analytic) stratified K\"ahler
polarization $P$ shortly.
For $\ell \geq s$, as a complex analytic space,
$N^{\roman{red}}$ comes down to
$\roman S_{\Bobb C}^2[\Bobb C^{\ell}]$
whereas, for
$\ell < s$, as a complex analytic space,
$N^{\roman{red}}$ may be described as a
complex determinantal variety
in
$\roman S_{\Bobb C}^2[\Bobb C^{\ell}]$, that is, as an affine
variety given by determinantal equations;
see e.~g. \cite\brunvett\ 
for determinantal varieties.
This may be deduced from standard geometric invariant theory
results combined with the standard description of the invariants
of the classical groups which, in turn, may be found e.~g. in
\cite\weylbook.
As a stratified symplectic space,
the singularity structure of 
$N^{\roman{red}}$ is finer
than that of the complex analytic structure, though:  
Once $\ell$ is fixed,
for $ s= \ell $,
the 
smooth structure $C^{\infty}(N_{\ell})$ and hence the
Poisson structure on
$N_{\ell} =N^{\roman{red}} \cong \Bobb C^d, d=\frac{\ell(\ell+1)}{2}$,
is not standard and,
as a stratified symplectic space,
$N_{\ell}$ has $\ell+1$ strata. For $s<\ell$,
the space $N^{\roman{red}} = N_s$ (say) 
may be described as the closure of a stratum
in $N_{\ell}$;
moreover, 
a system of $\ell$ particles 
in $\Bobb R^{s-1}$
being viewed as
a system of $\ell$ particles 
in $\Bobb R^s$
via the standard inclusion of
$\Bobb R^{s-1}$ into $\Bobb R^s$
yields an injection of $N_{s-1}$ into
$N_s$. Thus we get a sequence
$$
\{o\} =N_0 \subseteq N_1 \subseteq \dots N_{s-1} \subseteq N_s
\subseteq \dots \subseteq N_{\ell}
$$
of injections of normal (complex analytic
stratified) K\"ahler spaces
in such a way that,
for $1 \leq s \leq \ell$,
$N_{s-1} \subseteq N_s$ is the singular locus
of $N_s$ in the sense of stratified symplectic spaces,
and 
the stratified K\"ahler 
structure
on
$N_s$,
in particular the requisite Poisson structure,
is then simply obtained by restriction from $N_{\ell}$.
For example, for $\ell=1$,
$(N_1,C^{\infty}(N_1),\{\cdot,\cdot\})$
is just the reduced space and reduced Poisson algebra
of a system of a single particle in $\Bobb R^n$ ($n \geq 2$)
with angular momentum zero
explained in the Example 4.2 above.
For 
$\ell=s=2$,
the space
$N_2= N^{\roman{red}}$
is 
complex analytically
a copy of $\Bobb C^3$
which, as a stratified symplectic space, sits inside
$\fra{sp}(2,\Bobb R)$,
and we need ten generators to describe the Poisson
structure on $N_2$.
The reduced space              
$N_1$ for $\ell=2, s = 1$
is here complex analytically realized inside $N_2\cong \Bobb C^3$
as the quadric $Y^2=XZ$.
\smallskip
To introduce coordinates, and to spell out a description
of the complex analytic stratified K\"ahler polarizations, consider
the complexification
$\fra{sp}(\ell,\Bobb C)$
of
$\fra{sp}(\ell,\Bobb R)$; this complexification
sits inside the complex polynomial
algebra
$\Bobb C[z_1,\dots,z_{\ell},
\overline z_1,\dots,\overline z_{\ell}]$
as its homogeneous quadratic constituent.
The complexification
$\fra k^{\Bobb C} \cong \fra{gl}(\ell,\Bobb C)$
of
the maximal compact subalgebra
$\fra k=\fra{u}(\ell)$
of
$\fra{sp}(\ell,\Bobb R)$
is the span of the
$z_j \overline z_k$'s
and, with reference to the Cartan decomposition
$\fra{sp}(\ell,\Bobb R)= \fra{u}(\ell)\oplus \fra p$
of 
$\fra{sp}(\ell,\Bobb R)$,
the constituents $\fra p^+$ and $\fra p^-$
of the decomposition
$\fra p^{\Bobb C} = \fra p^+ \oplus \fra p^-$
are the spans of the
$z_j  z_k$'s
and
the $\overline z_j \overline z_k$'s,
respectively;
this gives an explicit description of 
$\fra p^+$ and $\fra p^-$
as
$\roman S_{\Bobb C}^2[\Bobb C^{\ell}]$
and
$\overline{\roman S_{\Bobb C}^2[\Bobb C^{\ell}]}$,
respectively.
Furthermore,
$\fra k=\fra{u}(\ell)$
sits inside
$\fra{sp}(\ell,\Bobb C)$
as the
real span of the
$z_j \overline z_k + \overline z_j  z_k$'s
and
$i(z_j \overline z_k - \overline z_j  z_k)$'s,
and
$\fra p$
sits 
inside
$\fra{sp}(\ell,\Bobb C)$
as the real span
of the
$z_j z_k + \overline z_j  \overline z_k$'s
and
$i(z_j  z_k - \overline z_j \overline z_k)$'s;
the assignment to a real symmetric
$(2\ell \times 2\ell)$-matrix
of the corresponding
complex
symmetric
$(\ell \times \ell)$-matrix
is given by the association
$$
z_j z_k + \overline z_j  \overline z_k
\longmapsto
z_j z_k,
\quad
i(z_j  z_k - \overline z_j \overline z_k)
\longmapsto
iz_j  z_k.
$$
The summands
$\fra p^+$ and $\fra p^-$
are the irreducible
$\fra k^{\Bobb C}$-representations
in
$\fra{sp}(\ell,\Bobb C)$
complementary to
$\fra k^{\Bobb C}$.
\smallskip
The homogeneous quadratic polynomials in
the variables $z_1,\dots,z_{\ell},
\overline z_1,\dots,\overline z_{\ell}$
yield coordinates on
$\fra{sp}(\ell,\Bobb R)$
and hence, via restriction, on
$N^{\roman{red}}$,
that is,
the smooth structure $C^{\infty}(N^{\roman{red}},\Bobb C)$
may be described as the algebra of smooth functions in these variables,
subject to the relations coming from the embedding of
$N^{\roman{red}}$ into 
$\fra{sp}(\ell,\Bobb R)$.
Now, the differentials $d(z_jz_k)$ of the coordinate functions
$z_jz_k$ $(1 \leq j,k \leq \ell$) 
(that is, of those coordinate functions which do not involve any of the
$\overline z_j$'s)
generate the corresponding
complex analytic stratified K\"ahler polarization
$P \subseteq \Omega^1(N^{\roman{red}},\Bobb C)$
as an
$(\Bobb C,C^{\infty}(N^{\roman{red}},\Bobb C))$-Lie subalgebra
of 
$\Omega^1(N^{\roman{red}},\Bobb C)_{\{\cdot,\}}$.
\smallskip
In \cite\kaehler, we developed a theory
of holomorphic nilpotent orbits of hermitian Lie algebras
and established the fact that the (topological) closure of any 
holomorphic nilpotent orbit inherits a normal (complex analytic
stratified) K\"ahler structure. 
The  space
$N^{\roman{red}}$,
realized  as the closure of a holomorphic nilpotent orbit
in
$\fra{sp}(\ell,\Bobb R)$,
is a special case thereof.
\smallskip\noindent
{\smc Example 6.2.} 
A variant of the above example
arises from the constrained system
of $\ell$ harmonic oscillators 
in $\Bobb R^s$ 
with total angular momentum zero
and constant energy.
Its unreduced
phase space
$Q$ is a 
copy of complex projective space $\Bobb P^{s\ell-1} \Bobb C$
of complex dimension $s\ell-1$.
For $\ell \geq s$, as a complex analytic space,
$Q^{\roman{red}}$ coincides with the (complex)
projectivization
$\Bobb P\roman S_{\Bobb C}^2[\Bobb C^{\ell}]$
of
$\roman S_{\Bobb C}^2[\Bobb C^{\ell}]$
whereas 
for $\ell < s$, as a complex analytic space,
$Q^{\roman{red}}$ may be described as a
complex projective determinantal variety
in
$\Bobb P\roman S_{\Bobb C}^2[\Bobb C^{\ell}]$.
In fact, the determinantal equations mentioned in Example 6.1
above are homogeneous
and yield the requisite homogeneous
equations for the present case.
In the same vein as before,
we get a sequence
$$
Q_1 \subseteq \dots Q_{s-1} \subseteq Q_s
\subseteq \dots \subseteq Q_{\ell} \cong \Bobb P^d \Bobb C,\ 
d = \frac{\ell(\ell+1)}2 -1,
$$
of injections of compact normal (complex analytic
stratified) K\"ahler spaces
in such a way that,
for $2 \leq s \leq \ell$,
$Q_{s-1} \subseteq Q_s$ is the singular locus
of $Q_s$ in the sense of stratified symplectic spaces,
each $Q_s$ being the closure of a stratum in $Q_{\ell}$,
and 
the stratified K\"ahler 
structure
on
$Q_s$,
in particular the requisite Poisson structure,
is simply obtained by restriction from $Q_{\ell}$.
Complex analytically, each $Q_s$ is a projective variety.
Again,
the smooth structure $C^{\infty}(Q_{\ell})$
and hence the Poisson structure on
$Q_{\ell}\cong \Bobb P^d \Bobb C$ $(s=\ell)$
is not the standard one (which
arises from the Fubini-Study metric on complex projective space)
and, as a stratified symplectic space,
$Q_{\ell}$ has $\ell$ strata.
For example, for
$\ell=s=2$,
the space
$Q_2$
is 
complex analytically
a copy of $\Bobb P^2\Bobb C$,
and the corresponding reduced space 
$Q_1$ (for $\ell=2, s = 1$),
which is abstractly just complex projective 1-space,
sits complex analytically inside $Q_2\cong \Bobb P^2\Bobb C$
as the projective conic $Y^2=XZ$.
These spaces are particular cases of a systematic class of examples
of {\it exotic projective varieties\/}, introduced and explored
in our paper \cite\kaehler.

\smallskip\noindent
{\smc Remark 6.3.\/}
Given a Lie group $G$, a smooth hamiltonian $G$-space,
and a real $G$-invariant polarization, 
the question arises whether the statement
of Proposition 5.1 still holds for this real polarization.
When we try to identify, on the reduced level,
a  stratified version of 
such a polarization, we may run into the following difficulty, though:
Under the circumstances of the Example 6.1,
let $\ell =1$, and
consider the
 vertical polarization on $N=\roman T^*\Bobb R^n$.
This polarization integrates to the foliation---even
fibration---defined by the
projection map from
$\roman T^*\Bobb R^n$ to $\Bobb R^n$.
This foliation is clearly $\roman O(n,\Bobb R)$-invariant and,
in terms of the standard coordinates $\bold q=(q^1,\dots,q^n)$
on $\Bobb R^n$,
a leaf is given by the equation $\bold q = \bold q_0$ where $\bold q_0$ 
is a constant.
We will now write the ordinary scalar product of two vectors
$\bold x$ and $\bold y$ as $\bold x\bold y$. 
With these preparations out of the way,
under the present circumstances, the assignment to
$(\bold q,\bold p) \in \roman T^*\Bobb R^n$
of
$x_1 = \bold q\bold q - \bold p \bold p$ and $x_2 = 2 \bold q \bold p$
yields a map from
$\roman T^*\Bobb R^n$ to the plane $\Bobb R^2$
which induces an isomorphism
of stratified symplectic spaces
from the reduced space
$N^{\roman{red}}$ onto the exotic plane described in the Example 4.2.
In particular, the radius function $r$ is given by
$r = \bold q\bold q + \bold p \bold p$.
Now $2\bold q\bold q = x_1+r$ 
whence, under reduction,
the leaf 
$\bold q = \bold q_0$
passes to the subspace of the plane 
given by the equation 
$$
x_1+r = 2\bold q_0\bold q_0 = c\ \text{(say).}
$$
For 
$\bold q_0 \ne 0$, in the $(x_1,x_2)$-plane, this is just the parabola
$
x_2^2 + 2c x_1 = c^2
$
since $r^2 =  x_1^2 +x_2^2$
while,
for $\bold q_0 = 0$, it is the
non-positive half $x_1$-axis.
The reason for this degeneracy is that
the leaf $\bold q=0$ is not transverse to the momentum mapping $\mu$
in the sense that, whatever $\bold p\in \Bobb R^n$,
$\mu(0,\bold p)= 0$  while
$\roman{ker}(d\mu(0,\bold p))$ and the tangent
space of the leaf at $(0,\bold p)$ do {\it not\/}
together span the tangent space of $N$ at $(0,\bold p)$.
Thus the reduced space is still foliated,
but one leaf is singular;
however even 
the restriction of this foliation to
the top stratum still has a singular leaf,
the negative half $x_1$-axis.
A little thought reveals that
this implies that this foliation cannot result from a
stratified real polarization,
the notion of stratified real polarization being defined in the same
fashion as a stratified complex polarization,
except that, 
on each stratum, the polarization should come down to a real polarization.
As a side remark we mention that the 
assignment to a leaf of its intersection point with the
non-negative $x_1$-axis
identifies the 
space of leaves with the non-negative $x_1$-axis,
and the latter in fact coincides with the orbit space
$\Bobb R^n \big / \roman O(n,\Bobb R)$.
This description visualizes the exceptional role played by
the 
non-positive $x_1$-axis.
The distribution parallel to this foliation, though,
is given by the hamiltonian vector field of the
function $\bold q\bold q$ in 
$C^{\infty}(N^{\roman{red}})$;
it has the form
$$
\{\bold q\bold q,-\} = \frac 12 \{x_1 + r,-\}
= -x_2 \frac \partial{\partial x_1} + 
(x_1 + r) \frac \partial{\partial x_2}
$$
and in particular vanishes on the non-positive half $x_1$-axis.
The function $\bold q\bold q$ generates a maximal abelian Poisson subalgebra
of $\left(C^{\infty}(N^{\roman{red}}),\{\cdot,\cdot\}^{\roman{red}}\right)$.
This phenomenon is typical for cotangent bundles
with a hamiltonian action of a Lie group arising from an action
of that group on the base with more than a single orbit type.
{\sl Thus we see that the question whether a polarization other than a
K\"ahler polarization descends to a stratified polarization on the
reduced level leads to certain delicacies, and we do not know to what
extent we can interpret it merely as a descent problem.\/}
\smallskip
The question whether, under suitable circumstances
so that in particular the reduced space is still a smooth manifold,
a real polarization descends has been studied in \cite\gotaythr.

\medskip\noindent {\bf 7. Prequantization on spaces with singularities}

\smallskip\noindent
To develop
prequantization over stratified symplectic spaces and
to describe the behaviour of prequantization under reduction, 
in our paper \cite\qr,
we  introduced {\it stratified prequantum 
modules\/} over stratified symplectic spaces.
A stratified
prequantum 
module
is defined in terms of the appropriate Lie-Rinehart algebra and
determines what we call a {\it costratified prequantum space\/}
but the two notions, though closely related, should not be confused.
\smallskip
Let $N$ be a stratified symplectic space, and let
$(A,\{\cdot,\cdot\})$
be its stratified symplectic Poisson algebra; 
a special case would
be
the ordinary symplectic Poisson algebra of 
a smooth symplectic manifold.
Consider the extension (2.4) of Lie-Rinehart algebras.
Given an
$(A\otimes \Bobb C)$-module
$M$,
we refer to
an 
$(A,\overline L_{\{\cdot,\cdot\}})$-module
structure 
$
\chi
\colon 
\overline L_{\{\cdot,\cdot\}}
\longrightarrow
\roman{End}_{\Bobb R}(M)
$
on $M$  
as a 
{\it prequantum module structure for\/}
$(A,{\{\cdot,\cdot\}})$
provided 
(i) the values of $\chi$ lie in
$\roman{End}_{\Bobb C}(M)$,
that is to say, the operators $\chi(a,\alpha)$ are complex linear
transformations,
and (ii)
for every $a\in A$,
$\chi(a,0) = i\,a\,\roman{Id}_M$ \cite{\souriau,\,\qr}.
\smallskip
We recall from \cite\poiscoho\ that
the assignment to $a \in A$ of
$(a,da) \in
\overline L_{\{\cdot,\cdot\}}$
yields a morphism $\iota$ of real Lie algebras
from
$A$ to
$\overline L_{\{\cdot,\cdot\}}$;
this reduces the construction of Lie algebra representations
of 
the Lie algebra which underlies the Poisson algebra
$A$  to the construction
of representations of
$\overline L_{\{\cdot,\cdot\}}$.
Thus, for any prequantum module $(M,\chi)$,
the composite of $\iota$ with $-i \chi$
is a representation
$a \mapsto \widehat a$
of the $A$ underlying real Lie algebra
on $M$, viewed as a complex vector space,
by $\Bobb C$-linear operators 
so that the constants  in $A$ 
act by multiplication
and so that the Dirac condition holds,
even though $M$ does not necessarily carry a Hilbert space structure.
These operators are given by the formula
$$
\widehat a (x) = \frac 1 i \chi(0,da) (x) + ax,
\quad
a \in A,\ x \in M.
\tag7.1
$$
\smallskip
For illustration, consider  an ordinary quantizable
symplectic manifold $(N,\sigma)$, with ordinary {\it prequantum bundle\/} 
$\zeta \colon \Lambda \to N$, that is, $\zeta$ is a complex line bundle
with
a connection
$\nabla$ whose curvature equals $-i\sigma$;
the
assignments
$\chi_{\nabla}(a,0) = i\,a\,\roman{Id}_M$ $(a \in A)$ and
$\chi_{\nabla}(0,\alpha) = \nabla_{\pi^{\sharp}(\alpha)}$
$(\alpha \in \Omega^1(N)_{\{\cdot,\cdot\}})$
then yield 
a prequantum module structure 
$$
\chi_{\nabla} 
\colon
\overline L_{\{\cdot,\cdot\}} @>>>
\roman{End}_{\Bobb C}(M)
\subseteq 
\roman{End}_{\Bobb R}(M)
$$
for
$(A,\{\cdot,\cdot\})$.
(Here
$\pi^{\sharp}\colon \Omega^1(N) \to \roman{Vect}(N)$
refers to the adjoint of the 2-form $\pi$
induced by the symplectic Poisson structure, cf. Section 2.)
This is just
the ordinary prequantization
construction
in another guise.
\smallskip
As before, consider a general stratified symplectic space $N$,
with stratified symplectic Poisson algebra
$(C^{\infty}(N),\{\cdot,\cdot\})$.
For each stratum $Y$, let
$(C^{\infty}(Y),\{\cdot,\cdot\}^Y)$
be its symplectic Poisson structure, and let
$$
0
@>>>
C^{\infty}(Y)
@>>>
\overline L_{\{\cdot,\cdot\}^Y}
@>>>
\Omega^1(Y)_{\{\cdot,\cdot\}^Y}
@>>>
0
$$
be the corresponding extension (2.4) of Lie-Rinehart algebras.
As in (1.5) of \cite\kaehler,
we define a {\it stratified prequantum module\/}
for $N$ to consist of
\newline\noindent
--- a prequantum module $(M,\chi)$ for
$(C^{\infty}(N),\{\cdot,\cdot\})$, together with,
\newline\noindent
--- for each stratum $Y$, a prequantum module structure
$\chi_Y$ for
$(C^{\infty}(Y),\{\cdot,\cdot\}^Y)$
on $M_Y =  C^{\infty}(Y) \otimes_{C^{\infty}(N)} M$
in such a way that the canonical
linear 
map 
of complex vector spaces
from $M$ to $M_Y$
is a
morphism 
of prequantum modules
from $(M,\chi)$ to $(M_Y,\chi_Y)$.
\smallskip
Given a stratified prequantum module $(M,\chi)$  for $N$,
when $Y$ runs through the strata of $N$, 
we refer to the system 
$$
\left(M_{\overline Y},
\chi_{\overline Y}\colon \overline L_{\{\cdot,\cdot\}^{\overline Y}}
@>>> \roman{End}_{\Bobb R}(M_{\overline Y})\right)
$$
of prequantum modules, together with,
for every
pair of strata $Y,Y'$ such that
$Y' \subseteq \overline Y$,
the induced morphism
$$
\left(M_{\overline Y},\chi_{\overline Y}\right)
@>>>
\left(M_{\overline Y'},\chi_{\overline Y'}\right)
$$
of prequantum modules,
as a {\it costratified prequantum space\/}.
More formally: Consider the category $\Cal C_N$ whose objects are the strata
of $N$ and whose morphisms are the inclusions $Y' \subseteq \overline Y$.
We define a 
{\it costratified complex vector space\/} on $N$ to be a contravariant functor
from $\Cal C_N$ to the category of complex vector spaces,
and a
{\it costratified prequantum space\/} on $N$ to be a 
costratified complex vector space
together with a compatible system of prequantum module structures.
Thus a 
 stratified prequantum module  $(M,\chi)$ 
for 
$(N,C^{\infty}(N),\{\cdot,\cdot\})$
determines a costratified prequantum space on $N$; see
the (1.4) and (1.5) of \cite\qr\ for details. 

\proclaim{Theorem 7.2}
Given 
a symplectic manifold $N$ with a hamiltonian action of
a compact Lie group $G$,
a $G$-equivariant prequantum bundle
$\zeta$
descends to a stratified prequantum module
$(\chi^{\roman{red}},M^{\roman{red}})$
for the stratified symplectic space
$(N^{\roman{red}},C^{\infty}(N^{\roman{red}}),\{\cdot,\cdot\}^{\roman{red}})$.
\endproclaim

\demo{Proof} See Theorem 2.1 of \cite\qr.\qed
\enddemo

Thus, phrased in the language of prequantum modules,
the relationship between the unreduced and reduced prequantum object
may be interpreted as one of  descent.
\smallskip
In particular, consider
a complex analytic stratified K\"ahler 
space
\linebreak
$(N,C^{\infty}(N),\{\cdot,\cdot\},P)$
(cf. Section 5 above or Section 2 of \cite\kaehler),
and
let $(M,\chi)$ be a stratified prequantum module
for $(C^{\infty}(N),\{\cdot,\cdot\})$.
We refer to
$(M,\chi)$ as a {\it complex analytic\/} stratified prequantum module
provided 
$M$ is the space of ($C^{\infty}(N)$-) sections of a complex $V$-line bundle
$\zeta$ on $N$ in such a way that
$P$ endows $\zeta$ via $\chi$ with a complex analytic structure.
If this happens to be the case,
$M^P$ necessarily amounts to the space of global sections of
the sheaf of germs of holomorphic sections of $\zeta$.
See Section 3 of
\cite\kaehler.

\beginsection 8. K\"ahler quantization and reduction

Let $G$ be a compact Lie group,
let $(N, \sigma,\mu)$ be a hamiltonian $G$-space
of the kind as that in the circumstances of Proposition 5.1, and suppose that
$N$ is quantizable.
Thus $N$ is, in particular, a positive K\"ahler manifold
with a holomorphic $G^{\Bobb C}$-action whose restriction to $G$
preserves the K\"ahler structure and is hamiltonian.
Write $P$ for the corresponding K\"ahler polarization,
necessarily $G$-invariant,
viewed as a  
$(\Bobb C, C^{\infty}(N,\Bobb C))$-Lie subalgebra of the
$(\Bobb C, C^{\infty}(N,\Bobb C))$-Lie algebra
$\Omega^1(N,\Bobb C)_{\{\cdot,\cdot\}}$, and 
let $\zeta$ be a 
prequantum bundle.
Via its connection, it acquires a
holomorphic structure, and the connection is the unique hermitian
connection for a corresponding hermitian structure.
The momentum mapping induces, in particular,
an infinitesimal action of the Lie algebra $\fra g$ of $G$
on $\zeta$ preserving the connection and hermitian structure.
Suppose that  this action lifts to a
$G$-action  on $\zeta$
preserving the connection and lifting the $G$-action on $N$.
For connected $G$,
the assumption that the $G$-action lift to
one on $\zeta$ is
(well known to be)
redundant
(since the infinitesimal action
is essentially given by the momentum mapping)
and
it will suffice
to replace $G$ by an appropriate covering group if need be.
Prequantization turns the space of smooth sections of $\zeta$
into a
prequantum module
for the ordinary smooth symplectic Poisson algebra of $N$. We
write this prequantum module as
$M$; it
inherits a $G$-action preserving the polarization $P$.
Hence the quantum module $M^P$, that is, the space
$\Gamma(\zeta)$
of global holomorphic sections of $\zeta$,
is a complex representation space for $G$.
This quantum module 
is the corresponding {\it unreduced\/}
quantum state space, 
except that there is no
Hilbert space structure present yet,
and 
{\it reduction after quantization\/},
for the {\it quantum state spaces\/}, amounts to taking
the space
$(M^P)^G$
of $G$-invariant holomorphic sections.
\smallskip
The projection map
from the space
of smooth $G$-invariant sections of $\zeta$
to
$M^{\roman{red}}$
restricts to a linear map
$$
\rho
\colon
\Gamma(\zeta)^G
@>>>
(M^{\roman{red}})^{P^{\roman{red}}}
\tag8.1
$$
of complex vector spaces, defined on
the space
$(M^P)^G=\Gamma(\zeta)^G$ of 
$G$-invariant holomorphic sections of $\zeta$.
Here and below
$P^{\roman{red}}$ 
refers to the
stratified K\"ahler polarization
the existence of which is asserted
 in Proposition 5.1 above,
$M^{\roman{red}}$
to the prequantum module for the
stratified K\"ahler space
mentioned in Theorem 7.2 (without having been made explicit there),
and
$(M^{\roman{red}})^{P^{\roman{red}}}$
to the
$P^{\roman{red}}$-invariants;
notice that  $P^{\roman{red}}$ is, in particular,
a Lie algebra whence it makes sense to talk about
$P^{\roman{red}}$-invariants.
A module of the kind
$(M^{\roman{red}})^{P^{\roman{red}}}$
is referred to as a {\it reduced quantum module\/}
in \cite\qr.
It acquires a
costratified Hilbert space structure, the requisite scalar products 
being induced
from appropriate hermitian structures via integration.
\smallskip
As far as the comparison of 
$G$-invariant unreduced quantum observables
and reduced quantum observables is concerned,
the statement that
{\it K\"ahler quantization commutes with reduction\/} amounts to the following,
cf. Theorem 3.6 in \cite\qr.

\proclaim{Theorem 8.2}
The data $(N,\sigma,\mu,M,P)$ being fixed so that,
in particular, $(N,\sigma,\mu)$ is a smooth hamiltonian $G$-space
structure on a quantizable positive K\"ahler manifold $N$
with a holomorphic $G^{\Bobb C}$-action,
the K\"ahler polarization being written as $P$,
let $f$ be a smooth $G$-invariant
function on $N$
which is quantizable in the sense that it preserves $P$.
Then its class $[f] \in C^{\infty}(N^{\roman{red}})
(=(C^{\infty}(N))^G /I^G)$
is quantizable,  
i.~e. preserves $P^{\roman{red}}$ and,
for every
$h \in (M^P)^G$,
$$
\rho(\widehat f (h))
=
\widehat {[f]} (\rho (h)).
\tag8.2.1
$$
\endproclaim

So far, we did not make any claim to the effect that the 
reduced quantum module
$(M^{\roman{red}})^{P^{\roman{red}}}$
amounts to a space of global holomorphic sections. 
We now recall that, under the circumstances of
Theorem 8.2, the
momentum mapping
is said to be {\it admissible\/}
provided, for every $m \in N$, the path of steepest descent through $m$
is contained in a compact set
\cite\sjamatwo, \cite\kirwaboo\ (\S 9).
For example, when the momentum mapping is proper it is admissible.
Likewise, the momentum mapping for a unitary representation
of a compact Lie group is admissible in this sense, see Example 2.1 in
\cite\sjamatwo.

\smallskip
The statement
\lq\lq {\it K\"ahler quantization commutes with reduction\/}\rq\rq\ 
is then completed by the following
two observations, cf.
\cite\qr\ ((3.7) and (3.8)).

\proclaim{Proposition 8.3} Under the circumstances of
Theorem {\rm 8.2},
when $\mu$ is admissible and when $N^{\roman{red}}$ has a top stratum
(i.~e. an open dense stratum),
for example when $\mu$ is proper,
the reduced stratified prequantum module $(M^{\roman{red}},\chi^{\roman{red}})$ 
is complex analytic, that is,
as a complex vector space, $M^{\roman{red}}$
amounts to the space of global holomorphic sections of a suitable
holomorphic $V$-line bundle on $N^{\roman{red}}$.
\endproclaim

The relevant $V$-line bundle on $N^{\roman{red}}$
may be found in \cite\sjamatwo\
(Proposition 2.11).

\proclaim{Theorem 8.4}{\rm [\sjamatwo] (Theorem 2.15)}
Under the circumstances of Theorem {\rm 8.2},
when the momentum mapping $\mu$
is proper, 
in particular, when $N$ is compact,
the map $\rho$ is an isomorphism of complex vector spaces.
\endproclaim

In this result, the properness condition, while sufficient,
is not necessary, that is, the map $\rho$ may be an isomorphism
without the momentum mapping being proper.
\smallskip
A version of 
Theorem 8.4 has been established in (4.15) of \cite\naramtwo;
cf. also \cite\sjamafou\ and the literature there, as well as
\cite\ramadthr\ and \cite\ctelethr\ 
for generalizations to higher dimensional sheaf cohomology.
\smallskip
\noindent
{\smc Remark 8.5.}
The statements of Theorems 8.2 and 8.4 are logically independent;
in particular the statement of Theorem 8.2
makes sense whether or not
$\rho$ is an isomorphism, and its proof does
not rely on $\rho$ being an isomorphism.
\smallskip
Thus we have consistent K\"ahler quantizations
on the unreduced and reduced spaces,
including a satisfactory treatment of observables,
as indicated by the formula (8.2.1). 
We have already pointed out in the introduction that examples 
in finite dimensions 
abound.
We hope 
that this kind of approach, suitably adapted,
will eventually yield
the quantization of certain
infinite dimensional systems
arising from field theory.
\smallskip
\noindent
{\smc Remark 8.6.}
{\smc Kempf's\/} descent lemma  
\cite\kempfone\ mentioned earlier
characterizes, among the
holomorphic $V$-line bundles which arise on
a geometric invariant theory quotient by the standard geometric invariant
theory construction,
those which are ordinary (holomorphic) line bundles.
In the 
circumstances of Theorem 8.4, complex analytically, the space
$N^{\roman{red}}$ 
is a geometric invariant theory quotient, and
the $V$-line bundle
which underlies the reduced quantum module
arises by the standard geometric invariant
theory construction.
Here the term \lq\lq descent\rq\rq\ is used in its strict sense.

\smallskip
\noindent
{\smc Illustration 8.7.}
Under the circumstances of the Example 6.2,
let $\Cal O(1)$ be the ordinary hyperplane bundle
on $Q=\Bobb P^{s\ell-1} \Bobb C$ and,
as usual,
for $k \geq 1$, write its $k$'th power as $\Cal O(k)$.
The 
unitary group $\roman U(s\ell)$ acts on
$\Bobb P^{s\ell-1} \Bobb C$ in a hamiltonian fashion having as
 momentum mapping 
 the familiar embedding of $\Bobb P^{s\ell-1} \Bobb C$
 into $\fra u(s\ell)^*$, and  the adjoint thereof yields a morphism of Lie
 algebras
 from $\fra u(s\ell)$ to
 $C^{\infty}(\Bobb P^{s\ell-1} \Bobb C)$, the latter being
 endowed with its symplectic Poisson structure
 coming from the Fubini-Study metric.
 It is a standard fact that,
 for $k \geq 1$, K\"ahler quantization, with reference
 to $k\omega$
 and $\Cal O(k)$ (where $\omega$ is the Fubini-Study form),
 yields the $k$'th symmetric power  of the
 standard representation  defining the 
 Lie algebra $ \fra u(s\ell)$, and this representation integrates to
 the $k$'th symmetric power  $E_s^k$ of the
 standard representation $E_s$ defining the
 group $\roman U(s\ell)$.
 (We use the subscript $-_s$ since here and below $\ell$ is fixed
 while $s$ varies.)
 The symmetry group $H=\roman O(s,\Bobb R)$ of the constrained system
 in (6.1) above
 appears as a subgroup of $\roman U(s\ell)$
 in an obvious fashion and,
 viewed as this subgroup, 
  $H$ centralizes the subgroup 
 $\roman U(\ell)=\roman {Sp}(\ell,\Bobb R)\cap \roman U(s\ell)$
 (the maximal compact subgroup $\roman U(\ell)$
 of $\roman {Sp}(\ell,\Bobb R)$); hence, for $k \geq 1$, the subspace
  $(E_s^k)^{H}$
   of $H$-invariants
  is a $\roman U(\ell)$-representation.
  On the other hand,
  with an abuse of notation,
  let $\Cal O(1)$ be the hyperplane bundle
on the reduced space
$Q_{\ell}=\Bobb P^d \Bobb C$, $d= \frac{\ell(\ell+1)}2-1$ and, 
for $k \geq 1$, let $\Cal O(k)$ be
its $k$'th power. The space of holomorphic sections thereof,
$\Gamma(\Cal O(k))$, amounts to the $k$'th symmetric power
$\roman S_{\Bobb C}^k[\fra p^*]$  of the dual of 
$\fra p=\roman S_{\Bobb C}^2[\Bobb C^\ell]$
(the space of homogeneous degree $k$ polynomial functions on $\fra p$).
For $1\leq s\leq \ell$ and $k \geq 1$,
maintain the notation  
$\Cal O(k)$  for the restriction of the $k$'th power of the hyperplane
bundle to $Q_s \subseteq Q_{\ell}$; for $s<\ell$,
the space of holomorphic sections   
$\widetilde E^k_s$ of
$\Cal O(k)$ is now a certain quotient of
$\widetilde E^k_{\ell}=\roman S_{\Bobb C}^k[\fra p^*]$ 
which will be made precise below.
\smallskip
For $1\leq s\leq \ell$, the composite of the embedding of
$N_s$ into $\fra{sp}(\ell,\Bobb R)^*$ with the
surjection from 
$\fra{sp}(\ell,\Bobb R)^*$ to $\fra{u}(\ell)^*$ induced from the
injection of $\fra{u}(\ell)$ into $\fra{sp}(\ell,\Bobb R)$
yields a map from
$N_s$ to $\fra{u}(\ell)^*$
which descends to 
a map from
$Q_s$ to $\fra{u}(\ell)^*$,
the adjoint of which induces a
morphism of Lie
 algebras
 from $\fra u(\ell)$ to
 $C^{\infty}(Q_s)$, the latter being 
 endowed with its {\it stratified symplectic Poisson structure\/}
 explained earlier.
 For $k \geq 1$, the space of sections $M^{\roman{red}}$ (cf. Theorem 7.2) 
 of $\Cal O(k)$, with reference to a
 $C^{\infty}(Q_s)$-module structure constructed in \cite\qr\ 
 and not made precise here,
 inherits a stratified prequantum module structure;
 and stratified K\"ahler quantization yields a $\roman U(\ell)$-representation
 on the space 
 $\widetilde E^k_s$, which amounts to that
 written earlier as $(M^{\roman{red}})^{P^{\roman{red}}}$, cf.
 (8.1), 
 in such a way that the map $\rho$ given as (8.1) above
 identifies the representation written above as $(E_s^{2k})^{H}$
 with $\widetilde E^k_s$;
 moreover, the spaces $(E_s^{2k-1})^{H}$ are zero.
 \smallskip
 We conclude with an explicit description of the spaces $(E_s^{2k})^{H}$
 or, equivalently, of the spaces $\widetilde E^k_s$:
 Introduce coordinates $x_1,\dots,x_\ell$ on $\Bobb C^\ell$. 
 These give rise to coordinates
$\{x_{i,j} = x_{j,i}; \, 1 \leq i,j \leq \ell\}$ on 
$\fra p=S_{\Bobb C}^2 [\Bobb C^{\ell}]$,
and the determinants
$$
\delta_1 = x_{1,1},
\ 
\delta_2 = \left | \matrix x_{1,1}& x_{1,2}\\
                           x_{1,2}& x_{2,2}
                   \endmatrix \right|,
\
\delta_3 = \left | \matrix x_{1,1}& x_{1,2} & x_{1,3}\\
                           x_{1,2}& x_{2,2} & x_{2,3}\\
                           x_{1,3}& x_{2,3} & x_{3,3}\\
                   \endmatrix \right|,
\
\text{etc.}
$$
are highest weight vectors for certain $\roman U(\ell)$-representations.
For $1 \leq s \leq r$ and $k \geq 1$,
the 
$\roman U(\ell)$-representation $\widetilde E^k_s$
is the sum of 
the irreducible representations having as
highest weight vectors 
the monomials
$$
\delta_1^{\alpha} \delta_2^{\beta} \ldots \delta_s^{\gamma},
\quad 
\alpha +2 \beta + \dots + s\gamma = k,
$$
and the morphism from $\widetilde E^k_s$
to
$\widetilde E^k_{s-1}$
coming from restriction 
from $Q_s$ to $Q_{s-1}$
is an isomorphism on 
the span of
those irreducible representations which do not  involve
$\delta_s$
and has the
span of
the remaining ones as its kernel.
In particular, this explains how $\widetilde E^k_s$ arises from
$\widetilde E^k_{\ell}=\roman S_{\Bobb C}^k[\fra p^*]$.
For $1 \leq s \leq \ell$, the system
$(\widetilde E^k_1,\widetilde E^k_2,\dots,\widetilde E^k_s)$
is an example of a costratified quantum space.

\smallskip
The alerted reader is invited to consult \cite\qr\ for more details.

\bigskip
\centerline{\smc References}
\medskip
\widestnumber\key{9999}
\ref \no  \abramars
\by R. Abraham and J. E. Marsden
\book Foundations of Mechanics
\publ Benjamin/Cum-
\linebreak
mings Publishing Company
\yr 1978
\endref

\ref \no \almekump
\by R. Almeida and A. Kumpera
\paper Structure produit dans la cat\'egorie des alg\'ebro\"\i des de Lie
\jour An. Acad. Brasil. Cienc.
\vol 53 
\yr 1981
\pages 247--250
\endref

\ref \no  \armcusgo
\by J. M. Arms,  R. Cushman, and M. J. Gotay
\paper  A universal reduction procedure for Hamiltonian group actions
\paperinfo in: The geometry of Hamiltonian systems, T. Ratiu, ed.
\jour MSRI Publ. 
\vol 20
\pages 33--51
\yr 1991
\publ Springer-Verlag
\publaddr Berlin $\cdot$ Heidelberg $\cdot$ New York $\cdot$ Tokyo
\endref

\ref \no \armamonc
\by J. M. Arms, J. E. Marsden, and V. Moncrief
\paper  Symmetry and bifurcation of moment mappings
\jour Comm. Math. Phys.
\vol 78
\yr 1981
\pages  455--478
\endref

\ref \no \atiy
\by M. F. Atiyah
\paper Complex analytic connections in fibre bundles
\jour Trans. Amer. Math. Soc.
\vol 85
\yr 1957
\pages 181--207
\endref

\ref \no \beilsche
\by A. A. Beilinson and V. V. Schechtmann
\paper Determinant bundles and Virasoro algebras
\jour Comm. Math. Physics 
\vol 118
\yr 1988
\pages 651--701
\endref

\ref \no \bkoucone
\by R. Bkouche
\paper Structures $(K,A)$-lin\'eaires
\jour C. R. A. S. Paris S\'erie A
\vol 262
\yr 1966
\pages 373--376
\endref

\ref \no \brunvett
\by W. Bruns and U. Vetter
\book Determinantal Rings
\bookinfo Lecture Notes in Mathematics, Vol. 1327
\yr 1988
\publ Springer-Verlag
\publaddr Berlin $\cdot$ Heidelberg $\cdot$ New York
\endref

\ref \no \canawein
\by Ana Canas da Silva and Alan Weinstein
\book Geometric models for Noncommutative Algebras
\bookinfo Berkeley Mathematical Lecture Notes, Vol. 10
\publ AMS
\publaddr Boston Ma
\yr 1999
\endref

\ref \no \chasetwo
\by S.~U. Chase
\paper Group scheme actions by inner automorphisms
\jour Comm. Alg.
\vol 4
\yr 1976
\pages 403--434
\endref

\ref \no \chemltwo
\by S. Chemla
\paper A duality property for complex Lie algebroids
\jour Math. Z.
\vol 232
\yr 1999
\pages 367--388
\endref

\ref \no \cheveile
\by C. Chevalley and S. Eilenberg
\paper Cohomology theory of Lie groups and Lie algebras
\jour  Trans. Amer. Math. Soc.
\vol 63
\yr 1948
\pages 85--124
\endref

\ref \no \craifern
\by M. Crainic and R. L. Fernandes
\paper Integrability of Lie brackets 
\jour Ann. of Math. (to appear)
\finalinfo{\tt math.DG/0105033}
\endref

\ref \no  \debarone
\by C.~M. de Barros
\paper Espaces infinit\'esimaux
\jour Cahiers Topologie G\'eom. diff\'erentielle
\vol 7
\yr 1964AA
\endref

\ref \no \diracone
\by P. A. M. Dirac
\book Lectures on Quantum Mechanics
\publ Belfer Graduate School of Science
\publaddr Yeshiva University, New York
\yr 1964
\endref

\ref \no \diractwo
\by P. A. M. Dirac
\paper Generalized Hamiltonian systems
\jour Can. J. of Math.
\vol 12
\yr 1950
\pages 129--148
\endref

\ref \no \dreznara
\by J.-M. Drezet and M.S. Narasimhan
\paper Groupe de Picard des vari\'et\'es de 
modules de fibr\'es semistables sur les courbes alg\'ebriques 
\jour Invent. Math. 
\vol 97 
\pages 53--94 
\yr 1989
\endref
 
\ref \no \fahim
\by A. Fahim
\paper Extensions galoisiennes d'alg\`ebres diff\'erentielles 
\jour Pacific J. Math. 
\vol 180 
\pages 7-40 
\yr 1997
\endref

\ref \no \fuchsste
\by B. Fuchssteiner
\paper The Lie algebra structure of degenerate Hamiltonian
and bi-hamiltonian systems
\jour Progr. Theor. Phys.
\vol 68
\yr 1982
\pages 1082--1104
\endref

\ref \no \gormacon
\by M. Goresky and R. MacPherson
\paper Intersection homology theory
\jour Topology
\vol 19
\yr 1980
\pages 135--162
\endref

\ref \no \gotaythr
\by M. J. Gotay
\paper Constraints, reduction, and quantization
\jour J. of Math. Phys.
\vol 27
\yr 1986
\pages 2051--2066
\endref

\ref \no  \guistetw
\by V. W. Guillemin and S. Sternberg
\paper Geometric quantization and multiplicities of group representations
\jour Invent. Math.
\vol 67
\yr 1982
\pages 515--538
\endref

\ref \no \hermaone
\by R. Hermann
\paper Analytic continuations of group representations. {\rm IV.}
\jour Comm. Math. Phys.
\vol 5
\yr 1967
\pages 131--156
\endref

\ref \no \herzone
\by J. Herz
\paper Pseudo-alg\`ebres de Lie
\jour C. R. Acad. Sci. Paris 
\vol 236
\yr 1953
\pages 1935--1937
\endref

\ref \no \higgmack
\by P. J. Higgins and K. Mackenzie
\paper Algebraic constructions in the category of Lie algebroids
\jour J. of Algebra
\vol 129 
\yr 1990
\pages 194--230
\endref

\ref \no \hinschon
\by V. Hinich and V. Schechtman
\paper Deformation theory and Lie algebra homology. {\rm I.}
\paperinfo {\tt alg-geom/9405013}
\jour Alg. Colloquium
\vol 4:2
\yr 1997
\pages 213--240
\moreref
II.
\pages 291--316
\endref

\ref \no \hochsfiv
\by G. Hochschild
\paper Simple algebras with purely inseparable splitting field of exponent 1
\jour  Trans. Amer. Math. Soc.
\vol 79
\yr 1955
\pages 477--489
\endref

\ref \no \hochkoro
\by G. Hochschild, B. Kostant, and A. Rosenberg
\paper Differential forms on regular affine algebras
\jour  Trans. Amer. Math. Soc.
\vol 102
\yr 1962
\pages 383--408
\endref

\ref \no \poiscoho
\by J. Huebschmann
\paper Poisson cohomology and quantization
\jour 
J. f\"ur die reine und angewandte Mathematik
\vol  408 
\yr 1990
\pages 57--113
\endref

\ref \no  \souriau
\by J. Huebschmann
\paper On the quantization of Poisson algebras
\paperinfo Symplectic Geometry and Mathematical Physics,
Actes du colloque en l'honneur de Jean-Marie Souriau,
P. Donato, C. Duval, J. Elhadad, G.M. Tuynman, eds.;
Progress in Mathematics, Vol. 99
\publ Birkh\"auser-Verlag
\publaddr Boston $\cdot$ Basel $\cdot$ Berlin
\yr 1991
\pages 204--233
\endref

\ref \no \srni
\by J. Huebschmann
\paper
Poisson geometry of certain
moduli spaces
\paperinfo
Lectures delivered at the \lq\lq 14th Winter School\rq\rq, Srni,
Czeque Republic,
January 1994
\jour Rendiconti del Circolo Matematico di Palermo, Serie II
\vol 39
\yr 1996
\pages 15--35
\endref

\ref \no \claustha
\by J. Huebschmann
\paper On the Poisson geometry of certain moduli spaces
\paperinfo in: Proceedings of an international workshop on
\lq\lq Lie theory and its applications in physics\rq\rq,
Clausthal, 1995,
H. D. Doebner, V. K. Dobrev, J. Hilgert, eds.
\publ World Scientific
\publaddr Singapore $\cdot$
New Jersey $\cdot$
London $\cdot$
Hong Kong 
\pages 89--101
\yr 1996
\endref

\ref \no \bv
\by J. Huebschmann
\paper Lie-Rinehart algebras, Gerstenhaber algebras, and Batalin-
Vilkovisky algebras
\jour Annales de l'Institut Fourier
\vol 48
\yr 1998
\pages 425--440
\linebreak
\finalinfo{\tt math.DG/9704005}
\endref
\ref \no \extensta
\by J. Huebschmann
\paper Extensions of Lie-Rinehart algebras and the Chern-Weil construction
\paperinfo Festschrift in honour of Jim Stasheff's 60'th anniversary
\jour Cont. Math.  
\vol 227 
\yr 1999 
\pages 145--176
\finalinfo{\tt math.DG/9706002}
\endref

\ref \no \duality
\by J. Huebschmann
\paper 
Duality for Lie-Rinehart algebras and the modular class
\jour Journal f\"ur die reine und angew. Math.
\vol 510
\yr 1999
\pages 103--159
\finalinfo{\tt math.DG/9702008}
\endref

\ref \no \banach
\by J. Huebschmann
\paper  Differential Batalin-Vilkovisky algebras arising from twilled 
Lie-Rinehart algebras
\jour Banach center publications
\vol 51
\yr 2000
\pages 87--102
\endref 

\ref  \no  \oberwork
\by J. Huebschmann
\paper Singularities and Poisson geometry of certain representation spaces
\paperinfo in: Quantization of Singular Symplectic Quotients,
N. P. Landsman, M. Pflaum, M. Schlichenmaier, eds.,
Workshop, Oberwolfach,
August 1999,
Progress in Mathematics, Vol. 198
\publ Birkh\"auser-Verlag
\publaddr Boston $\cdot$ Basel $\cdot$ Berlin
\yr 2001
\pages 119--135
\finalinfo{\tt math.DG/0012184}
\endref

\ref \no \twilled
\by J. Huebschmann
\paper Twilled Lie-Rinehart algebras and differential Batalin-Vilkovisky 
algebras
\paperinfo {\tt math.DG/9811069}
\endref

\ref  \no \kaehler
\by J. Huebschmann
\paper K\"ahler spaces, nilpotent orbits, and singular reduction
\linebreak
\paperinfo {\tt math.DG/0104213}
\endref

\ref \no \severi
\by J. Huebschmann
\paper Severi varieties and holomorphic nilpotent orbits
\linebreak
\paperinfo {\tt math.DG/0206143}
\endref

\ref \no \qr
\by J. Huebschmann
\paper K\"ahler reduction and quantization
\paperinfo {\tt math.SG/0207166}
\endref

\ref \no \illusboo
\by L. Illusie
\book Complexe cotangent et d\'eformations. {\rm II}
\bookinfo Lecture Notes in Mathematics
 No. 238
\publ Springer-Verlag
\publaddr Berlin $\cdot$  Heidelberg $\cdot$ New York
\yr 1972
\endref

\ref \no \jacobfiv
\by N. Jacobson
\paper An extension of Galois theory to non-normal and non-separable fields
\jour Amer. J. Math.
\vol 66
\yr 1944
\pages 1--29
\endref 

\ref \no \kambtond
\by F. W. Kamber and Ph. Tondeur
\paper Invariant differential operators and the cohomology
of Lie algebra sheaves
\jour Memoirs of the Amer. Math. Soc.
\vol 113
\yr 1971
\publ Amer. Math. Soc.
\publaddr Providence, R. I
\pages  
\endref

\ref \no \kasstora
\by D. Kastler and R. Stora
\paper A differential geometric setting for BRS transformations 
and anomalies. {\rm I.II.}
\jour J. Geom. Phys.
\vol 3
\yr 1986
\pages  437--482; 483--506
\endref

\ref \no \kempfone
\by G. Kempf
\paper Instability in invariant theory 
\jour Ann. of Math.
\vol 108
\yr 1978
\pages 299--316
\endref 

\ref \no \kemneone
\by G. Kempf and L. Ness
\paper The length of vectors in representation spaces
\jour Lecture Notes in Mathematics
\vol 732
\yr 1978
\pages 233--244
\paperinfo Algebraic geometry, Copenhagen, 1978
\publ Springer-Verlag
\publaddr Berlin $\cdot$ Heidelberg $\cdot$ New York
\endref

\ref \no \kirilone
\by A. A. Kirillov
\paper Unitary representations of nilpotent Lie groups
\jour Uspehi Mat. Nauk.
\vol 17
\yr 1962
\pages 57--101
\moreref
\jour
Russ. Math. Surveys
\vol 17
\yr 1962
\pages 57--101
\endref

\ref \no \kirwaboo
\by F. Kirwan
\book Cohomology of quotients in symplectic and algebraic geometry
\publ Princeton University Press
\publaddr Princeton, New Jersey
\yr 1984
\endref

\ref \no \kosmagri
\by Y. Kosmann-Schwarzbach and F. Magri 
\paper Poisson-Nijenhuis structures
\jour  Annales Inst. H. Poincar\'e S\'erie A (Physique th\'eorique)
\vol 53
\yr 1989
\pages 35--81
\endref

\ref \no \kostaone
\by B. Kostant
\paper Quantization and unitary representations
\jour Lecture Notes in Math.
\vol 170
\yr 1970
\pages 87--207
\paperinfo In:
Lectures in Modern Analysis and Applications, III, ed. C. T. Taam
\publ Springer-Verlag
\publaddr Berlin $\cdot$ Heidelberg $\cdot$ New York
\endref

\ref \no \kostetwo
\by B. Kostant and S. Sternberg
\paper Anti-Poisson algebras and current algebras
\paperinfo unpublished manuscript
\yr 1990
\endref

\ref \no \kralyvin
\by I. S. Krasil'shchik, V. V. Lychagin, and A. M. Vinogradov
\book Geometry of Jet Spaces and Nonlinear Partial Differential Equations
\bookinfo Advanced Studies in Contemporary Mathematics, vol. 1
\publ Gordon and Breach Science Publishers
\publaddr New York, London, Paris, Montreux, Tokyo
\yr 1986
\endref

\ref \no \lermonsj
\by E. Lerman, R. Montgomery and R. Sjamaar
\paper Examples of singular reduction
\paperinfo Symplectic Geometry,
Warwick, 1990,  D. A. Salamon, editor, 
London Math. Soc. Lecture Note 
Series
\vol 192
\yr 1993
\pages  127--155
\publ Cambridge University Press
\publaddr Cambridge, UK
\endref

\ref \no \lichn
\by A. Lichnerowicz
\paper Les vari\'et\'es de Poisson et leurs alg\`ebres de Lie
associ\'ees
\jour J. Diff. Geo.
\vol 12
\yr 1977
\pages 253--300
\endref

\ref \no \mackbook
\by K. Mackenzie
\book Lie groupoids and Lie algebroids in differential geometry
\bookinfo London Math. Soc. Lecture Note Series, vol. 124
\publ Cambridge University Press
\publaddr Cambridge, England
\yr 1987
\endref

\ref \no \mack
\by K. Mackenzie
\paper Lie algebroids and Lie pseudoalgebras
\jour Bull. London Math. Soc. 
\vol 27 (2)
\pages 97 -- 147
\yr 1995 
\endref

\ref \no \marswein
\by J. Marsden and A. Weinstein
\paper Reduction of symplectic manifolds with symmetries
\jour Rep. on Math. Phys.
\vol 5
\yr 1974
\pages 121--130
\endref

\ref \no \masspete
\by W. S. Massey and F. P. Petersen
\paper The cohomology structure of certain fibre spaces.I
\jour Topology
\vol 4
\yr 1965
\pages  47--65
\endref

\ref \no \naramtwo
\by M. S. Narasimhan and T. R. Ramadas
\paper Factorization of generalized theta functions
\jour Inventiones
\vol 114
\yr 1993
\pages 565-623 
\endref

\ref \no \nelsoboo
\by E. Nelson
\book Tensor Analysis
\publ Princeton University Press
\publaddr Princeton, N. J.
\yr 1967
\endref

\ref \no \palaione
\by R. S. Palais
\paper The cohomology of Lie rings
\jour  Proc. Symp. Pure Math.
\vol III
\yr 1961
\pages 130--137
\paperinfo Amer. Math. Soc., Providence, R. I.
\endref

\ref \no \prad
\by J. Pradines
\paper Th\'eorie de Lie pour les groupo\"\i des diff\'erentiables.
Relations entre propri\'et\'es locales et globales 
\jour  C. R. Acad. Sci. Paris S\'erie A
\vol 263
\yr 1966
\pages 907--910
\endref

\ref \no \praditwo
\by J. Pradines
\paper Th\'eorie de Lie pour les groupo\"\i des diff\'erentiables.
Calcul diff\'erentiel dans la cat\'egorie des
groupo\"\i des  infinit\'esimaux
\jour C. R. Acad. Sci. Paris S\'erie A
\vol 264
\yr 1967
\pages 245--248
\endref

\ref \no \ramadthr
\by T. R. Ramadas
\paper Factorization of generalised theta functions {\rm II:}
The Verlinde formula
\jour Topology
\vol 35
\yr 1996
\pages  641--654
\endref 

\ref \no \rinehone
\by G. Rinehart
\paper Differential forms for general commutative algebras
\jour  Trans. Amer. Math. Soc.
\vol 108
\yr 1963
\pages 195--222
\endref

\ref \no \segal
\by I. E. Segal
\paper Quantization of non-linear systems
\jour J. of Math. Phys.
\vol 1
\yr 1960
\pages 468--488
\endref

\ref \no \sjamatwo
\by R. Sjamaar
\paper Holomorphic slices, symplectic reduction, and multiplicities of 
representations
\jour Ann. of Math.
\vol 141
\yr 1995
\pages 87--129
\endref

\ref \no \sjamafou
\by R. Sjamaar
\paper Symplectic reduction and Riemann-Roch formulas for multiplicities
\jour Bull. Amer. Math. Soc.
\vol 33
\yr 1996
\pages 327--338
\endref

\ref \no \sjamlerm
\by R. Sjamaar and E. Lerman
\paper Stratified symplectic spaces and reduction
\jour Ann. of Math.
\vol 134
\yr 1991
\pages 375--422
\endref

\ref \no \sniabook
\by J. \'Sniatycki 
\book Geometric quantization and quantum mechanics
\bookinfo Applied Mathematical Sciences
 No.~30
\publ Springer-Verlag
\publaddr Berlin $\cdot$ Heidelberg $\cdot$ New York
\yr 1980
\endref

\ref \no \sniatone
\by J. ~\'Sniatycki 
\paper Constraints and quantization
\paperinfo in: Nonlinear partial differential operators
and quantization procedures,
Clausthal 1981, eds. S.~I.~ Anderson and H.~D.~Doebner
\jour Lecture Notes in Mathematics, No.~1037
\pages 301--334
\publ Springer-Verlag
\publaddr Berlin $\cdot$ Heidelberg $\cdot$ New York
\yr 1983
\endref

\ref \no \sniawein
\by J. ~\'Sniatycki and A. Weinstein
\paper Reduction and quantization for singular moment mappings
\jour Lett. Math. Phys.
\vol 7
\yr 1983
\pages 155--161
\endref

\ref \no \sourione
\by J. M. Souriau
\paper 
Quantification g\'eom\'etrique
\jour Comm. Math. Physics
\vol 1
\yr 1966
\pages 374--398
\endref

\ref  \no \stashfiv
\by J. D. Stasheff
\paper Homological reduction of constrained Poisson algebras
\jour J. of Diff. Geom. 
\vol 45
\yr 1997
\pages 221--240
\endref

\ref \no \nteleman
\by N. Teleman
\paper A characteristic ring of a Lie algebra extension
\jour Accad. Naz. Lincei. Rend. Cl. Sci. Fis. Mat. Natur. (8)
\vol 52
\yr 1972
\pages 498--506 and 708--711
\endref

\ref  \no \ctelethr
\by C. Teleman
\paper The quantization conjecture revisited
\jour Ann. of Math.
\vol 152
\yr 2000
\pages 1--43
\endref

\ref \no \weylbook
\by H. Weyl
\book The classical groups
\publ Princeton University  Press
\publaddr Princeton, New Jersey
\yr 1946
\endref

\ref \no \woodhous
\by N. M. J. Woodhouse
\book Geometric quantization
\bookinfo Second edition
\publ Clarendon Press
\publaddr Oxford
\yr 1991
\endref

\enddocument